\documentclass[times,final]{elsarticle}

\usepackage{jcomp}
\usepackage{framed,multirow}

\usepackage{graphicx}
\usepackage{subfigure}

\usepackage{amsmath,amssymb}
\usepackage{amsthm}

\usepackage{mathtools}
\usepackage{enumitem}

\usepackage{algorithm}
\usepackage{algorithmicx}
\usepackage{algpseudocode}

\usepackage{hyperref}
\usepackage{float}
\usepackage{epsfig}
\usepackage{epstopdf}
\usepackage{array,booktabs}
\usepackage{tabularx}
\usepackage[pagewise]{lineno}  
\usepackage{MnSymbol}
\providecommand{\U}[1]{\protect\rule{.1in}{.1in}}

\usepackage{longtable}

\newtheorem{lemma}{Lemma}
\newtheorem{corollary}{Corollary}

\newtheorem{assumption}{Assumption}
\newtheorem{remark}{Remark}
\newtheorem{definition}{Definition}%
\newtheorem{theorem}{Theorem}%
\newtheorem{proposition}{Proposition}%

\usepackage{MnSymbol}

\usepackage{url}
\usepackage{xcolor}
\definecolor{newcolor}{rgb}{.8,.349,.1}

\journal{}

\begin{document}

\verso{Dmitrii Chaikovskii, Ye Zhang}

\begin{frontmatter}

\title{Convergence analysis for forward and inverse problems in singularly perturbed time-dependent reaction-advection-diffusion equations\tnoteref{tnote1} }%
\tnotetext[tnote1]{Please cite to this paper as published in: \\  Journal of Computational Physics, 470:111609, 2022, \\ https://doi.org/10.1016/j.jcp.2022.111609.}

\author[1]{Dmitrii \snm{Chaikovskii}}
 \ead{mitichya@yandex.ru}
\author[2,1]{Ye \snm{Zhang}\corref{cor1}}
\ead{ye.zhang@smbu.edu.cn}
\cortext[cor1]{Corresponding author.   }

\address[1]{Shenzhen MSU-BIT University, Shenzhen, 518172, China}
\address[2]{School of Mathematics and Statistics, Beijing Institute of Technology, Beijing, 100081, China}


\begin{abstract}
In this paper, by employing the asymptotic expansion method, we prove the existence and uniqueness of a smoothing solution for a time-dependent nonlinear singularly perturbed partial differential equation (PDE) with a small-scale parameter. As a by-product, we obtain an approximate smooth solution, constructed from a sequence of reduced stationary PDEs with vanished high-order derivative terms. We prove that the accuracy of the constructed approximate solution can be in any order of this small-scale parameter in the whole domain, except a negligible transition layer. Furthermore, based on a simpler link equation between this approximate solution and the source function, we propose an efficient algorithm, called the asymptotic expansion regularization (AER), for solving nonlinear inverse source problems governed by the original PDE. The convergence-rate results of AER are proven, and the a posteriori error estimation of AER is also studied under some a priori assumptions of source functions. Various numerical examples are provided to demonstrate the efficiency of our new approach.
\end{abstract}


\end{frontmatter}



\section{Introduction}
\label{sec:Introduction}

This paper is mainly concerned with the usage of the asymptotic expansion method and the a posteriori error estimation for a nonlinear inverse source problem in time-dependent reaction–diffusion–advection equations. To introduce the underlying idea, we consider the following one-dimensional problem as an example:\\
(\textbf{IP}): Given noisy data $\{u^\delta(x,t_0), w^\delta(x,t_0)\}$ of $\{u(x,t), \frac{\partial u}{\partial x}(x,t)\}$ at the $n$ location points $\{x_1, \cdots, x_n\}$ and at only one time point $t_0$, find the source function $f(x)$ such that $(u,f)$ satisfies the nonlinear autowave model
\begin{align} \label{equat1}
\begin{cases}
\displaystyle \mu \frac{\partial^{2}u}{\partial x^{2}}-\frac{\partial u}{\partial t}=-ku \frac{\partial u}{\partial x}+ f(x), \quad x\in(0,1)\equiv\Omega, \quad  t\in(0, T],\\
\displaystyle u(0, t)=u^{l}, \quad  u(1, t)=u^{r}, \quad   t\in[0, T]\equiv \bar{\mathcal{T}},\\
u(x,0)=u_{init}(x), \quad  x\in[0,1]=\bar{\Omega},
\end{cases}
\end{align}
where $0<\mu \ll 1$ is the turbulent-diffusion coefficient, $u(x,t)$ is the dimensionless pollutant density, $x$ is the spatial coordinate, $t$ is the time variable, $k$ is the positive coefficient of distribution of a pollutant in the environment, and $f(x)$ is a function representing the intensity and location of the pollutant source. For simplicity, the left and right boundary conditions $u^{l}$ and $u^{r}$ are assumed to be constants. Model \eqref{equat1} assumes that the propagation speed does not depend on the water flow rate but depends only on the amount of the pollutant, and that pollutant dissipates quickly (for example, the pollutant could be noise pollution of the water or the spread of electricity in the water). In this work, the focus is on the speed, location, and width of the border between two regions -- a region with a high concentration of a pollutant and a region with an acceptable concentration. We assume that, at the middle point of this border (named the \emph{transition point}), which defines the location of the border, the concentration of the pollutant is equal to the critical value of the pollutant density in the medium.

Reaction–diffusion–advection PDEs with small parameters arise in a wide range of scientific disciplines, such as astrophysics \cite{BerrymanHolland1978}, biology  \cite{PattersonWagner2012,DoOwida2011,BodnarSequeira2008,HidalgoTello2014}, liquid chromatography  \cite{zhang2016regularization,zhang2016,LinZhang2018,ChengLin2018}, and industrial and environmental problems \cite{PhysRevE.70.026307,2004evapormetals,BerrymanHolland1978,2014evolutiondisp}. Specific examples of interest here are the reaction–diffusion–advection models for predicting the spread of environmental pollution, as shown in \cite{b01,b02}. However, as we are investigating rapidly dispersing pollutants in this paper, we suggest using an autowave approach with the model \eqref{equat1}. The application of the autowave model to the reaction–diffusion problem was considered in \cite{b04}, the authors proposing a model to predict the growth of the city of Shanghai in subsequent years. The authors also created a similar model that allowed them to analyze, over time, the displacement of the boundary between the two regions -- the region with a high population density and the region with a low population density.

Numerical schemes for solving the forward reaction–diffusion–advection equations comprise (a) finite difference methods \cite{AnguelovLubuma2003,ClaveroGracia2005,Mickens2000}, (b) finite element methods \cite{ArayaBehrens2005,FrancaValentin2000,IdelsohnNigro1996}, and (c) finite volume methods \cite{TitarevToro2002}. Asymptotic methods are especially attractive for partial differential equations with small parameters, e.g. the model \eqref{equat1} with small diffusion parameter $\mu$, since this technique allows us to find approximate solutions of singularly perturbed boundary-value problems and express these solutions in terms of known functions or quadratures from them, and also allows us to prove the existence and uniqueness of these solutions \cite{Tikhonov1948,ButuzovVasileva1970,AntipovLevashova2018}.
 In particular, the closer the small parameter is to zero, the more effective the asymptotic methods are, as the system becomes very difficult for the traditional numerical solution to handle. Another advantage of asymptotic methods is that the numerical solution is pointwise, while the asymptotic solution is smooth. Hence, the first objective of this paper is to develop an appropriate asymptotic method for the PDE model \eqref{equat1}.

The practical aim of the nonlinear inverse problem (\textbf{IP}) is to help eliminate negative environmental impacts. Because of the increasing amount of cargo being transported by water and the corresponding increase in the number of cargo ships, the problem of noise pollution of water is a very pressing one. The diesel engines and propellers of cargo ships generate high noise levels \cite{b1,b2}, and this noise pollution significantly increases the levels of low-frequency ambient noise. Even marine invertebrates such as crabs are affected by ship noise \cite{b4}, and noise pollution could have killed some species of whales that came ashore after exposure to the loud sounds of military sonar \cite{b5}. It should be noted that these kinds of inverse problems arise in many physical and engineering problems, such as \cite{math9040342,w12092415,math8050777}. Inverse source problems for other control equations can be found, for example, in \cite{Yamamoto1995}. In particular, inverse problems for the Burgers-type equations were recently considered in \cite{NefedovVolkov2020}, where the authors restored the modular-type source, and \cite{shishlenin2020}, where the authors reconstructed the initial condition $u_{init}(x)$ from the observation data of the transition layer. It is well known that the inverse source problem is severely ill-posed because of the unboundedness of its inversion operator (see \cite{Isakov1990} for a detailed discussion of the theoretical aspects of this problem). Therefore, for the problem of noisy boundary data, regularization methods should be employed to obtain meaningful source functions. With Tikhonov regularization, the inverse source problem (\textbf{IP}) can be converted to the following minimization problem:
\begin{equation} \label{LSM0}
\min_{f\in Q_{ad}} \sum^n_{i=0} \left\{ \left[u(x_i,t_0) - u^\delta(x_i,t_0) \right]^2 + \left[\frac{\partial u}{\partial x}(x_i,t_0) - w^\delta (x_i,t_0) \right]^2  \right\} + \varepsilon \mathcal{R}(f),
\end{equation}
where $u$ solves \eqref{equat1} with the given $f$, $\mathcal{R}(f)$ denotes the regularization term, and $\varepsilon>0$ is the regularization parameter. $Q_{ad}$ is an admissible set, incorporating the a priori information about the source function.

It is well known that the regularization parameter plays a crucial role in solving an ill-posed problem. In practice, in order to select an optimal value for the regularization parameter, one needs to repeatedly solve the forward problem \eqref{equat1}, which is usually time-consuming, especially for large-scale problems. The main idea in this paper is to use the asymptotic analysis to reduce the original nonlinear singularly perturbed problem to simpler problems without small parameters and high-order derivative terms while obtaining a sufficiently accurate solution. It should be noted that a similar idea was used in \cite{shishlenin2018}, where the coefficient inverse problem was considered for a nonlinear singularly perturbed reaction–diffusion–advection equation. Later, by using the same methodology, the authors in \cite{Paper-7} reconstructed the boundary condition from the observation of the reaction front for the similar PDE model. However, all the works mentioned here (i.e. \cite{Paper-1,Paper-6,shishlenin2020,shishlenin2018,Paper-7}) focus on the numerical implementation of the method. Hence, the main aim of this work is to establish a rigorous mathematical theory, i.e. convergence results as well as the (computable) a posteriori error estimation, for asymptotic-analysis-based inversion approaches. Note that, besides the presented model \eqref{equat1}, the framework proposed in this paper can also be applied to various linear and nonlinear inverse problems in singularly perturbed PDEs, e.g. inverse source problems in parabolic or hyperbolic singularly perturbed PDEs \cite{Atifi2018}, parameter-identification problems in singularly perturbed PDEs \cite{VolkovNevedofCoefInvProbl2020,Paper-1,Paper-6,2015Mustonen}, etc.

The remainder of the paper is structured as follows. Section \ref{statementresults} states the main results of the research, including the approximation results for both forward and inverse problems of PDE \eqref{equat1}. Section \ref{derivationAndProofs} covers the construction of asymptotic solutions and the technical proofs of the main theoretical results. Section \ref{simulation} presents some numerical experiments using our method, for both forward and inverse problems. Finally, Section \ref{Conclusion} concludes the paper.

\section{Statement of main results}\label{statementresults}

Table \ref{NotationTable} below lists the notations and abbreviations that will be used in this section.

\small{
\begin{longtable}{|l|l|l|}
\hline \multicolumn{1}{|c|}{\textbf{Notation}} & \multicolumn{1}{c|}{\textbf{Description }} & \multicolumn{1}{c|}{\textbf{Reference}} \\ \hline

\hline
\hspace{-3.5mm} \begin{tabular}{l}  $\mu$   \end{tabular} & \hspace{-3.5mm} \begin{tabular}{l}  Small parameter in the model, $ 0<\mu \ll 1 $       \end{tabular}   &  Eq.\eqref{equat1}  \\
 \hline
\renewcommand{\arraystretch}{1.2} \hspace{-4.5mm}  \begin{tabular}{l} $\Omega$, $\bar{\mathcal{T}}$ \end{tabular} & \hspace{-3.5mm} \begin{tabular}{l} Spacial and time domains: \\ $\Omega=(0,1)$, $\bar{\mathcal{T}}=[0, T]$ \end{tabular} & Eq.\eqref{equat1} \\
 \hline
\hspace{-3.5mm} \begin{tabular}{l}    $u^{l,r}$ \end{tabular} & \hspace{-3.5mm} \begin{tabular}{l}  Left and right boundary conditions \end{tabular}  & Eq.\eqref{equat1} \\
\hline
\renewcommand{\arraystretch}{1.2}  
\hspace{-4.5mm} \begin{tabular}{l}  $w^\delta(x,t)$   \end{tabular} & \hspace{-3.5mm} \begin{tabular}{l}  Noisy data of $\frac{\partial u}{\partial x}  (x,t) $       \end{tabular}   &  Eq.\eqref{noisyData1}  \\
   \hline
\hspace{-3.5mm} \begin{tabular}{l}  $x_{t.p.}$ \end{tabular} & \hspace{-3.5mm} \begin{tabular}{l}  Point at the middle of the transition   layer \end{tabular}  & Eq.\eqref{expantionxtp} \\
 \hline
\hspace{-3.5mm} \begin{tabular}{l}  $v_{t.p.}$ \end{tabular} & \hspace{-3.5mm} \begin{tabular}{l} Velocity at point $x_{t.p.}$, i.e.   $v_{t.p.}= \frac{d x_{t.p.}}{dt}$  \end{tabular} & Eq.\eqref{expantionv} \\
 \hline
\hspace{-3.5mm} \begin{tabular}{l}  Superscript  ``\,$^{l,r}$\,'' \end{tabular} & \hspace{-3.5mm} \begin{tabular}{l} Describes functions on the left and right, 
respectively, \\ relative to the point $x_{t.p.}$
\end{tabular} & Eq.\eqref{u}   \\
   \hline
\hspace{-3.5mm} \begin{tabular}{l} Subscript ``\,$_{0,1, \cdots}$\,'' \end{tabular} & \hspace{-3.5mm} \begin{tabular}{l} Describes the order of approximation of    the asymptotic solution \end{tabular} & Eq.\eqref{expantionxtp}-\eqref{transitionexpansions} \\
 \hline
\hspace{-4.5mm} \renewcommand{\arraystretch}{1.2}  \begin{tabular}{l}   $\bar{\Omega}^{l}$, $\bar{\Omega}^{r}$ \end{tabular} & \hspace{-3.5mm} \begin{tabular}{l}  Left and right regions relative to point $x_{t.p.}$: \\ $\bar{\Omega}^{l}=[0,x_{t.p.}(t,\mu)]$, $\bar{\Omega}^{r}=[x_{t.p.}(t,\mu),1]$
\end{tabular}   & Eq.\eqref{u}  \\
 \hline
\hspace{-3.5mm} \begin{tabular}{l} $X_{n}(t,\mu)$ \end{tabular} & \hspace{-3.5mm} \begin{tabular}{l} Approximation terms of $x_{t.p.}$  up to  order $n$ \end{tabular} & Eq.\eqref{XtpOfOrderN}  \\
 \hline
\hspace{-3.5mm} \begin{tabular}{l}  $U_{n}$ \end{tabular} & \hspace{-3.5mm} \begin{tabular}{l}  The $n$-order asymptotic solution \end{tabular} & Eq.\eqref{asymptoticnorder} \\
 \hline
\hspace{-3.5mm} \begin{tabular}{l} $\bar{u}^{l,r}$ \end{tabular} & \hspace{-3.5mm} \begin{tabular}{l} Regular functions describing the   solution away from the point $x_{t.p.}$ \end{tabular} & Eq.\eqref{ulr}  \\
 \hline
\hspace{-3.5mm} \begin{tabular}{l} $\varphi^{l,r}$ \end{tabular} & \hspace{-3.5mm} \begin{tabular}{l} Zero-order approximation of the   regular functions \end{tabular} & Eq.\eqref{eq7},\eqref{u0regu}  \\
 \hline
\hspace{-3.5mm} \begin{tabular}{l} $Q^{l,r}$ \end{tabular} &  \hspace{-3.5mm} \begin{tabular}{l} Transition-layer functions describing  the solution \\ near the point $x_{t.p.}$ \end{tabular} & Eq.\eqref{ulr} \\
 \hline
\hspace{-3.5mm}  \begin{tabular}{l}  $ \Delta x$  \end{tabular} & \hspace{-3.5mm} \begin{tabular}{l}   The width of the transition layer,  $ \Delta x \sim \mu \vert \ln \mu \vert $  \end{tabular}   &  Eq.\eqref{transitionVariables} \\
 \hline
\hspace{-3.5mm} \begin{tabular}{l} $ \xi$   \end{tabular} & \hspace{-3.5mm} \begin{tabular}{l}  Extended variable, $\xi = \frac{x-x_{t.p.}(t,\mu)}{\mu}$       \end{tabular}   & Eq.\eqref{xi} \\
 \hline
  \hspace{-4.5mm}  \renewcommand{\arraystretch}{1.2}
  \begin{tabular}{l} $u^\varepsilon (x,t_0)$ \end{tabular}  & \hspace{-3.5mm} \begin{tabular}{l}     The smoothed approximate data at time $t_0$ \end{tabular}   & Eq.\eqref{uAlpha} \\
 \hline
  \hspace{-4.5mm}  \renewcommand{\arraystretch}{1.2}
  \begin{tabular}{l} $f^*(x) $ \end{tabular}  & \hspace{-3.5mm} \begin{tabular}{l}      The exact source function   \end{tabular}   & Eq.\eqref{f0Ineq} \\
 \hline
 \hspace{-4.5mm}  \renewcommand{\arraystretch}{1.2}
  \begin{tabular}{l}  $f^\delta(x)$ \end{tabular}  & \hspace{-3.5mm} \begin{tabular}{l}      The reconstructed approximate \\ source function \end{tabular}   & Eq.\eqref{uAlpha1} \\
  \hline
 \hspace{-4.5mm}  \renewcommand{\arraystretch}{1.2}
  \begin{tabular}{l}  $Q^{1,2}_{ad}$ \end{tabular}  & \hspace{-3.5mm} \begin{tabular}{l} The set of approximate source functions      \end{tabular}   & Eq.\eqref{SetApp} \\
    \hline
 \hspace{-4.5mm}  \renewcommand{\arraystretch}{1.2}
  \begin{tabular}{l}  $\Delta_1$  \end{tabular}  & \hspace{-3.5mm} \begin{tabular}{l}   The relative a posteriori error    \end{tabular}   & Eq.\eqref{posterioriErr2} \\
  \hline
 \hspace{-4.5mm}  \renewcommand{\arraystretch}{1.2}
  \begin{tabular}{l}  $\Delta_2(x)$  \end{tabular}  & \hspace{-3.5mm} \begin{tabular}{l}   The pointwise error estimate    \end{tabular}   & Eq.\eqref{pointwiseErr} \\
      \hline
 \hspace{-4.5mm}  \renewcommand{\arraystretch}{1.2}
  \begin{tabular}{l}  $f^{low}(x), f^{up}(x)$  \end{tabular}  & \hspace{-3.5mm} \begin{tabular}{l}   The lower and upper solutions   \end{tabular}   & Eq.\eqref{UpandLowSol} \\
   \hline
\caption{Notations and references to their definitions.}
\label{NotationTable}
\end{longtable}
}

First, we list the assumptions under which the asymptotic solution exists:

\begin{assumption}
\label{A1}
$u^{l}<0$, $u^{r}>0$ and $u^r-u^l > 2\mu^2 $.
\end{assumption}

\begin{assumption}
\label{A2}
\begin{align}\label{A2bounds1}
\frac{k}{2} (u^{r})^2 >  \int_0^1  \max (0,f(x )) d x,\quad \frac{k}{2} (u^{l})^2 > - \int_0^1  \min (0,f(x ))  d x.
\end{align}

Note that, if the source function $ f(x) $ is positive, Assumption \ref{A2} can be replaced with the boundedness of the $L^1$ energy of source function $f(x)$, i.e.
\begin{align} \label{A2bounds}
\|f\|_{L^1(\Omega)} \leq \frac{k}{2} (u^{r})^2.
\end{align}
\end{assumption}

\begin{assumption}
\label{A3}
For all $t\in \bar{\mathcal{T}}$, $0<x_0(t)<1$, where $x_0$ is the zero approximation of $x_{t.p.}$; see \eqref{expantionxtp}.
\end{assumption}

\begin{assumption}
\label{A4}
$u_{init}(x)=U_{n-1}(x,0)+\mathcal{O}(\mu^n)$, where the asymptotic solution $U_{n-1}$ will be constructed later; see, e.g., Theorem \ref{MainThm}.
\end{assumption}

Under Assumptions \ref{A1} and \ref{A2}, zero-order regular functions $\varphi^{l} (x)$ and $\varphi^{r}(x)$, which will be used for asymptotic construction (refer to \eqref{u0regu}), can be expressed explicitly from equation \eqref{LeftRightSolution} in the following form:
\begin{align}
\label{eq7}
\varphi^{l} (x) = -\sqrt{\frac{2}{k} \int_0^x f(s )  d s + (u^{l})^2},\quad
\varphi^{r} (x) = \sqrt{ (u^{r})^2 - \frac{2}{k} \int_x^1  f(s ) d s}.
\end{align}

Assumption \ref{A3} defines the location of the internal transition layer within the specified region $\Omega$, while Assumption \ref{A4} means that at $t=0$ the transition layer has already been formed, and the initial function $u_{init}$ already has the transition layer in the vicinity of the point $x_{0}^{*}:= x_0(0)$. It should be noted that for some special $f(x)$, for which equation \eqref{eq10} has an explicit solution, Assumption \ref{A3} can be replaced by the requirements for the coefficients and boundary conditions, which are more reasonable in practice. However, in this paper, since we are mainly interested in solving inverse problems (i.e., the reconstruction of $f(x)$), equation \eqref{eq10} cannot be analytically solved without the information of $f(x)$. Therefore, in the general case, Assumption \ref{A3} cannot be replaced by a simpler condition. Nevertheless, an empirical simplification of Assumption \ref{A3} is discussed in Remark \ref{remarkx0}.  

We are now in a position to provide the main result for the forward problem (the meanings and descriptions of some notations can be found in Table \ref{NotationTable} and Section \ref{constructionOfAsymptotic}, respectively).
\begin{theorem} \label{MainThm}
Suppose that $f(x) \in C^1(\bar{\Omega})$, $u_{init}(x) \in C^1 (\bar{\Omega})$, and $\mu\ll 1$. Then, under Assumptions \ref{A1}–\ref{A4}, the boundary-value problem \eqref{equat1} has a unique smooth solution with an internal transition layer. In addition, the $n$-order asymptotic solution $U_{n}$ has the following representation ($\xi_{n} = (x-X_{n}(t,\mu)) / \mu $):
\begin{align} \label{asymptoticnorder}
 U_{n}(x,t,\mu)= \begin{cases}
 \displaystyle U_{n}^{l}(x,t,\mu)=\sum_{i=0}^{n} \mu^{i} (\bar{u}_{i}^{l}(x)+Q_{i}^{l}(\xi_{i},t)), \quad (x,t) \in \bar{\Omega}^{l} \times \bar{\mathcal{T}} , \\
 \displaystyle  U_{n}^{r}(x,t,\mu)=\sum_{i=0}^{n} \mu^{i} (\bar{u}_{i}^{r}(x)+Q_{i}^{r}(\xi_{i},t)), \quad (x,t) \in \bar{\Omega}^{r}\times \bar{\mathcal{T}},
  \end{cases}
 \end{align}
which  gives an approximation of the solution of problem \eqref{equat1} uniformly on the interval $x\in \bar{\Omega}$. Moreover, the following asymptotic estimates hold:
\begin{equation} \label{NorderEstim1}
\forall (x,t)\in \bar{\Omega}  \times \bar{\mathcal{T}}:~ \lvert u(x,t)-U_{n}(x,t,\mu) \rvert \leq \mathcal{O} ( \mu^{n+1}) ,
\end{equation}
\begin{equation} \label{NorderEstim2}
\forall t \in \bar{\mathcal{T}}:~ \lvert x_{t.p.}(t,\mu)-X_{n}(t,\mu) \rvert \leq \mathcal{O} ( \mu^{n+1}),
\end{equation}
\begin{align}\label{NorderEstim3}
\forall (x,t) \in \bar{\Omega}\backslash \{ X_{n}(t,\mu) \}  \times \bar{\mathcal{T}}:~ \left \lvert \frac{\partial u(x,t)}{\partial x}-\frac{\partial  U_{n}(x,t,\mu) }{\partial x} \right \rvert \leq \mathcal{O} ( \mu^{n}).
\end{align}

\end{theorem}

\begin{corollary}
\label{Corollary1}
(Zeroth approximation) Under the assumptions of Theorem \ref{MainThm}, the zero-order asymptotic solution $U_{0}$ has the following representation:
 \begin{align} \label{eq001}
U_{0}(x,t)=\begin{cases}
\varphi^{l}(x)+Q_{0}^{l}(\xi_{0},t) , \quad \ (x,t)\in \bar{\Omega}^{l}  \times \bar{\mathcal{T}},\\
\varphi^{r}(x)+Q_{0}^{r}(\xi_{0},t) , \quad \  (x,t)\in \bar{\Omega}^{r}  \times \bar{\mathcal{T}},
\end{cases}
\end{align}
where $\xi_{0}=(x-x_{0}(t))/\mu$. Moreover, the following holds:
\begin{equation} \label{eq002}
\forall (x,t)\in \bar{\Omega}  \times \bar{\mathcal{T}}:~ \vert u(x,t)-U_{0}(x,t)\vert =\mathcal{O}(\mu) ,
\end{equation}
\begin{equation} \label{eq003}
\forall t \in \bar{\mathcal{T}}:~ \lvert x_{t.p.}(t,\mu)-x_{0}(t) \rvert =\mathcal{O}(\mu) .
\end{equation}

Furthermore, outside the narrow region $(x_{0}(t)- \Delta x /2, x_{0}(t)+\Delta x /2)$ with $\Delta x\sim\mu \lvert \ln\mu \rvert $, there exists a constant $C$ independent of $\mu$, $x$ and $t$ such that the following inequalities hold:
\begin{equation} \label{eq004}
\lvert u(x,t)-\varphi^{l}(x) \rvert \leq C \mu, \qquad  (x,t) \in [0, x_{0}(t)- \Delta x /2]  \times \bar{\mathcal{T}},
\end{equation}
\begin{equation} \label{eq005}
\lvert u(x,t)-\varphi^{r}(x) \rvert  \leq C \mu, \qquad  (x,t) \in [x_{0}(t)+\Delta x /2, 1]  \times \bar{\mathcal{T}},
\end{equation}

\begin{equation}\label{0orderEstim1}
\left \lvert \frac{\partial u(x,t)}{\partial x}-\frac{d  \varphi^{l}(x)}{d x} \right \rvert  \leq C \mu,  \quad  (x,t) \in [0, x_{0}(t)- \Delta x /2]   \times \bar{\mathcal{T}},
\end{equation}
\begin{equation}\label{0orderEstim2}
\left \lvert \frac{\partial u(x,t)}{\partial x}-\frac{d \varphi^{r}(x)}{d x} \right \rvert \leq C \mu,  \quad  (x,t) \in [x_{0}(t)+\Delta x /2, 1]  \times \bar{\mathcal{T}}.
\end{equation}
\end{corollary}

Corollary \ref{Corollary1} follows directly from Theorem \ref{MainThm}. The inequalities \eqref{0orderEstim1}–\eqref{0orderEstim2} in Corollary \ref{Corollary1} can be obtained by taking into account the fact that the transition-layer functions are decreasing functions with respect to $\xi_0$ and are sufficiently small at the boundaries of the narrow region $(x_0(t)-\Delta x/2, x_0(t)+\Delta x/2)$, i.e. equation \eqref{Qlr-mu}.
From Corollary \ref{Corollary1}, it follows that the solution can be approximated by regular functions of zero order everywhere, except for a thin transition layer.

\begin{remark} 
\label{Remark1}
The point-wise error $E_{n}(x,t,\mu)$ of the asymptotic approximation $U_{n} (x,t,\mu)$ is defined by
\begin{align*} 
E_{n}(x,t,\mu) = \left|u(x,t)-U_{n} (x,t,\mu)\right| = \begin{cases}
 \displaystyle E_{n}^{l}(x,t,\mu)=\left|u(x,t)-U_{n}^{l} (x,t,\mu)\right|, \quad (x,t) \in \bar{\Omega}^{l} \times \bar{\mathcal{T}} , \\
 \displaystyle  E_{n}^{r}(x,t,\mu)=\left|u(x,t)-U_{n}^{r} (x,t,\mu)\right|, \quad (x,t) \in \bar{\Omega}^{r}\times \bar{\mathcal{T}}.
  \end{cases}
\end{align*}

According to Theorem \ref{MainThm}, the constant in the error estimation of asymptotic solution, i.e. \eqref{NorderEstim1}, can be approximately (and numerically) estimated by 
\begin{align*} 
C = \sup_{(x,t)\in \bar{\Omega}  \times \bar{\mathcal{T}}} \frac{E_{n}(x,t,\mu)}{\mu^{n+1}}, \qquad \text{for~} 0<\mu \ll 1.
\end{align*}

If the asymptotic solution converges, the error is equal to the remainder of the asymptotic series after truncation:
\begin{align*} 
E_{n}^{l,r}(x,t,\mu) = \left| \sum_{i=n+1}^{\infty} \mu^{i} \left( \bar{u}_i^{l,r}(x)+Q_i^{l,r}(\xi_i,t) \right) \right|, \quad (x,t) \in \bar{\Omega}^{l,r} \times \bar{\mathcal{T}} .
\end{align*}

Of course, for a fixed $\mu$, the asymptotic series may diverge, then there exists an index $N(\mu) $ of the order of the asymptotic approximation at which the error of the asymptotic solution is minimal \cite{truncation2005}. The corresponding error value is called the optimal error $E_{\text{opt}}(x,t,\mu)= \min_{N}\left|u(x,t)-U_N (x,t,\mu)\right|$. Moreover, the point-wise error for the divergent asymptotic solution can be calculated through the Borel summation \cite[Note 4.96]{Costin2008AsymptoticsAB}:
\begin{align*} 
E_{n}^{l,r}(x,t,\mu) = \left|  \int_0^\infty \left( \sum_{i=0}^{\infty}  \frac{\bar{u}_{i}^{l,r}(x)+Q_{i}^{l,r}(\xi_{i},t) }{i!}(\mu z)^i \right)e^{-z} dz - U_n^{l,r} (x,t,\mu) \right|, \quad (x,t) \in \bar{\Omega}^{l,r} \times \bar{\mathcal{T}}.
\end{align*}

Thus, for the zeroth asymptotic approximation, the constant in the error estimate \eqref{eq004} for the convergent and divergent asymptotic solution, respectively, can be numerically calculated by:
\begin{align*} 
\sup_{(x,t)\in \bar{\Omega}^{l}  \times \bar{\mathcal{T}}} \left| \frac{ Q_0^{l}(\xi_i,t)+\sum_{i=1}^{\infty}  \mu^{i } \left( \bar{u}_i^{l}(x)+Q_i^{l}(\xi_i,t) \right)}{\mu}  \right| = C, \qquad \text{for~} 0<\mu \ll 1,
\end{align*}

\begin{align*} 
\sup_{(x,t)\in \bar{\Omega}^{l}  \times \bar{\mathcal{T}}} \left| \frac{1}{\mu} \left( \int_0^\infty \left( \sum_{i=0}^{\infty}  \frac{ \bar{u}_{i}^{l}(x)+Q_{i}^{l}(\xi_{i},t) }{i!}(\mu z)^i \right)e^{-z} dz - \varphi^{l}(x) \right) \right| = C, \qquad \text{for~} 0<\mu \ll 1.
\end{align*}

The equations defining the constant $C$ in the estimate \eqref{0orderEstim1} in the case of convergent and divergent $x$-derivative of the asymptotic solution have the following representation, respectively:

\begin{align*} 
\sup_{(x,t)\in [0, x_{0}(t)- \Delta x /2]   \times \bar{\mathcal{T}}} \left| \frac{ \frac{\partial }{\partial x} \left( Q_0^{l}(\xi_i,t)+\sum_{i=1}^{\infty} \mu^{i } \left( \bar{u}_i^{l}(x)+Q_i^{l}(\xi_i,t)\right)\right)}{\mu}  \right| = C, \qquad \text{for~} 0<\mu \ll 1,
\end{align*}

\begin{align*} 
\sup_{(x,t)\in [0, x_{0}(t)- \Delta x /2]   \times \bar{\mathcal{T}}} \left| \frac{1}{\mu} \left( \int_0^\infty \left( \sum_{i=0}^{\infty}  \frac{\frac{\partial }{\partial x} \left(\bar{u}_{i}^{l}(x)+Q_{i}^{l}(\xi_{i},t) \right) }{i!}(\mu z)^i \right)e^{-z} dz - \frac{d \varphi^{l}(x)}{d x}   \right) \right| = C, \qquad \text{for~} 0<\mu \ll 1.
\end{align*}

In the same way we obtain the constants $C$ from the equations \eqref{eq005} and \eqref{0orderEstim2}.

According to Remark \ref{Remark1}, the following lemma holds.
\begin{lemma} 
\label{LemmaC}
Suppose that for a.e. $t\in \bar{\mathcal{T}}$,
\begin{align*} 
& \left\|\bar{u}_1^{l}+Q_1^{l} \right\|_{L^p(\bar{\Omega}^{l})} 
\cdot \left\| \frac{\partial }{\partial x} (\bar{u}_1^{l}+Q_1^{l}) \right\|_{L^p(0,x_0 - \Delta x/2)} 
\cdot  \left\| \frac{\partial Q_0^{l} }{\partial x}  \right\|_{L^p(0,x_0 - \Delta x/2)} \\ & \quad \cdot \left\|\bar{u}_1^{r}+Q_1^{r} \right\|_{L^p(\bar{\Omega}^{r})} 
\cdot \left\| \frac{\partial }{\partial x} (\bar{u}_1^{r}+Q_1^{r}) \right\|_{L^p(x_0 + \Delta x/2, 1)} 
\cdot \left\| \frac{\partial Q_0^{r} }{\partial x} \right\|_{L^p(x_0 + \Delta x/2,1)}
\neq 0,
\end{align*}
and the asymptotic solution $U_n(x,t,\mu)$ and its $x$-derivative converge. Then, for small enough $\mu<1$ such that
\begin{align*} 
\left\|\sum_{i=2}^{\infty} \mu^{i-1} \left( \bar{u}_i^{l}+Q_i^{l} \right) \right\|_{L^p(\bar{\Omega}^{l})} \leq \left\|\bar{u}_1^{l}+Q_1^{l} \right\|_{L^p(\bar{\Omega}^{l})}, \quad 
\left\|\sum_{i=2}^{\infty} \mu^{i-1} \frac{\partial }{\partial x} (\bar{u}_i^{l}+Q_i^{l}) \right\|_{L^p(0,x_0 - \Delta x/2)} \leq \left\| \frac{\partial }{\partial x} (\bar{u}_1^{l}+Q_1^{l}) \right\|_{L^p(0,x_0 - \Delta x/2)},
\end{align*}
\begin{align*} 
\left\|\sum_{i=2}^{\infty} \mu^{i-1} \left( \bar{u}_i^{r}+Q_i^{r} \right) \right\|_{L^p(\bar{\Omega}^{r})} \leq \left\|\bar{u}_1^{r}+Q_1^{r} \right\|_{L^p(\bar{\Omega}^{r})}, \quad 
\left\|\sum_{i=2}^{\infty} \mu^{i-1} \frac{\partial }{\partial x} (\bar{u}_i^{r}+Q_i^{r}) \right\|_{L^p(x_0 + \Delta x/2,1)} \leq \left\| \frac{\partial }{\partial x} (\bar{u}_1^{r}+Q_1^{r}) \right\|_{L^p(x_0 + \Delta x/2,1)},\\
\mu \left\| \frac{\partial Q_1^{l} }{\partial x}  \right\|_{L^p(0,x_0 - \Delta x/2)} \leq \left\| \frac{\partial Q_0^{l} }{\partial x}  \right\|_{L^p(0,x_0 - \Delta x/2)}, \quad \mu \left\| \frac{\partial  Q_1^{r} }{\partial x} \right\|_{L^p(x_0 + \Delta x/2,1)} \leq \left\| \frac{\partial Q_0^{r} }{\partial x} \right\|_{L^p(x_0 + \Delta x/2,1)},
\end{align*}
the following estimate holds
\begin{equation*} 
\left\| \varphi^{l}(x) - u (x,t) \right\|_{W^{1,p}(0,x_0 - \Delta x/2)} + \left\| \varphi^{r}(x) - u (x,t) \right\|_{W^{1,p}(x_0 + \Delta x/2,1)} \leq C \mu,
\end{equation*}
for a.e. $t\in \bar{\mathcal{T}}$. Here,
\begin{multline*}
C = 2\left(1 + \left\|\bar{u}_1^{l}+Q_1^{l} \right\|_{L^p(\bar{\Omega}^{l})}  + \left\| \frac{\partial \bar{u}_{1}^{l}}{\partial x}   \right\|_{L^p(0,x_0 - \Delta x/2)}+ \left\|\bar{u}_1^{r}+Q_1^{r} \right\|_{L^p(\bar{\Omega}^{r})}  + \left\| \frac{\partial \bar{u}_{1}^{r}}{\partial x}  \right\|_{L^p(x_0 + \Delta x/2,1)} \right)\\ + 3k\left( \left\|P^{l}\right\|_{L^p(0,x_0 - \Delta x/2)} +  \left\|P^{r}\right\|_{L^p(x_0 + \Delta x/2,1)} \right).
\end{multline*}
\end{lemma}

\end{remark}

Based on the asymptotic approximation of the forward problem, we proceed to construct an efficient inversion algorithm for (\textbf{IP}).
The main idea behind our new inversion algorithm is to replace the original governing  reaction–diffusion–advection equation \eqref{equat1} with a simpler relation \eqref{LeftRightSolution}. To this end, let $u (x,t)$ be the solution of PDE \eqref{equat1}, and define the pre-approximate source function $f_0$ as
\begin{equation} \label{f0}
f_0(x) = k u (x,t_0) \frac{u (x,t_0)}{dx}.
\end{equation}

\begin{proposition} 
\label{ProAsympErr}
Let $f^*$ be the exact source function, satisfying the original governing equation \eqref{equat1}. Then, under Assumption \ref{A2} and assumptions of Lemma \ref{LemmaC},
there exists a constant $C_1$ such that
\begin{equation}
\label{f0Ineq}
\|f^* - f_0\|_{L^p(\Omega)} \leq C_1 \mu \lvert \ln \mu \rvert, \quad \forall  p\in(0,+\infty).
\end{equation}
\end{proposition}

Suppose we have the deterministic noise model
\begin{equation}
\label{noisyData1}
\lvert u(x_i,t_0)-u^{\delta}_i \rvert \leq \delta, \quad \left \lvert \frac{\partial u(x_i,t_0)}{\partial x}-w^{\delta}_i \right \rvert \leq \delta,
\end{equation}
between the noisy data $\{u^{\delta}_i, w^{\delta}_i\}$ and the exact data $\{u(x_i,t_0), \frac{\partial u}{\partial x}(x_i,t_0)\}$ at time $t_0$ and at grid points $\Theta := \{0=x_0 < x_1 < \cdots < x_n =1 \}$ with maximum mesh size $h:= \max\limits_{i\in \{0, \cdots,n-1\}} \{x_{i+1} - x_{i}\}$.

Now, we restore the source function $f^\delta(x)$ according to the least-squares problem:
\begin{equation}
\label{uAlpha1}
f^\delta(x) = \mathop{\arg\min}_{\begin{subarray}{c} f\in C^1(0,1) \end{subarray} } \frac{1}{n+1} \sum^{n}_{i=0} \left( f(x_i) - k u^{\delta}_i w^{\delta}_i \right)^2.
\end{equation}

For the case in which we have only the noisy measurement $u^\delta_i$, we replace the values $\{ u^{\delta}_i, w^{\delta}_i \}$ in \eqref{uAlpha1} with smoothed quantities $\{u^\varepsilon, \frac{\partial u^\varepsilon}{\partial x}\}$. The function $u^\varepsilon (x,t_0)$ is constructed according to the following optimization problem:
\begin{equation}
\label{uAlpha}
u^\varepsilon(x,t_0) = \mathop{\arg\min}_{\begin{subarray}{c} s\in C^1(0,1)\end{subarray} } \frac{1}{n+1} \sum^{n}_{i=0} \left( s(x_i,t_0)-u^\delta_i(t_0) \right)^2 \\ + \varepsilon(t_0) \left\| \frac{\partial^2 s(x,t_0)}{\partial x^2} \right\|^2_{L^2(\Omega)},
\end{equation}
where the regularization parameter $\varepsilon(t_0)$ is chosen according to the discrepancy principle, i.e. the minimizing element $u^\varepsilon(x,t_0)$ of \eqref{uAlpha} satisfies $ \displaystyle \frac{1}{n+1} \sum^{n}_{i=0} \left( u^\varepsilon(x_i,t_0)-u^\delta_i(t_0) \right)^2 = \delta^2$.

\begin{proposition} 
\label{NoisyErr}
Suppose that, for a.e. $t\in\bar{\mathcal{T}}$, $u(\cdot,t) \in C^2(\bar{\mathcal{T}}, L^2(\Omega))$. Let $u^\varepsilon(x,t)$ be the minimizer of problem \eqref{uAlpha}, with $t_0$ replaced with $t$. Then, for a.e. $t\in\bar{\mathcal{T}}$  ,
\begin{equation}
\label{NoisyErrIneq}
\|u^\varepsilon(\cdot,t) - u(\cdot,t)\|_{H^1(0,1)} \leq 10\sqrt{2} \left( h  \left\| \frac{\partial^2 u(x,t)}{\partial x^2} \right\|_{L^2(\Omega)} + \sqrt{\delta} \left\|\frac{\partial^2 u(x,t)}{\partial x^2} \right\|^{1/2}_{L^2(\Omega)} \right). 
\end{equation}
\end{proposition}

\begin{remark}
Since the norms in the right hand side of \eqref{NoisyErrIneq} have large values in the transition layer, we use asymptotic analysis to exclude the transition layer from the optimization problem \eqref{uAlpha} and perform the minimization separately on the left and right regions $[0,x_0(t_0)-\Delta x /2 ]$ and $[x_0(t_0)+\Delta x/2 ,1]$. 
\end{remark}

With the help of Propositions \ref{ProAsympErr} and \ref{NoisyErr} and standard argument in the approximation theory, it is not difficult to construct the following theorem:

\begin{theorem} \label{ErrSource}
$f^\delta$, defined in \eqref{uAlpha1}, is a stable approximation of the exact source function $f^*$ for problem (\textbf{IP}). Moreover, it has the convergence rate
\begin{equation}
\label{ErrSourceIneq0}
\|f^* - f^\delta\|_{L^2(\Omega)} = \mathcal{O} ( \mu \lvert \ln \mu \rvert + h + \sqrt{\delta} ).
\end{equation}
Furthermore, if $\mu=\mathcal{O}(\delta^{\epsilon + 1/2})$ ($\epsilon$ is any positive number) and $h=\mathcal{O}(\sqrt{\delta})$, the following holds:
\begin{equation*}
\|f^* - f^\delta\|_{L^2(\Omega)} = \mathcal{O}(\sqrt{\delta}) .
\end{equation*}
\end{theorem}

At the end of this section, we study the error estimates for the obtained approximate source function $f^\delta$  with a priori information, i.e. $f^*\in Q_{ad}$, where the admissible set $Q_{ad}$ is assumed to be some compact sets. In this paper, we focus on the set of monotonic functions $M_1$ and the set of (piece-wise) convex functions $M_2$. It can be clearly shown that both $M_1$ and $M_2$ are compact sets in $L^2(\Omega)$ (see \cite{tihgonstya,TitarenkoLeonovYagola2002,TitarenkoYagola2008}).

First, define the set of approximate source functions as
\begin{equation}
\label{SetApp}
Q^{1,2}_{ad} := \left\{ f\in M_{1,2} \cap M_0  :~\mid f(x_i) -  k u^\delta_i w^\delta_i \mid \leq C_u \delta +  C_1 \mu \lvert \ln \mu \rvert, ~ i=0, \cdots, n \right\},
\end{equation}
where $C_u:= 2(\|u\|_{L^\infty(\Omega\times \mathcal{T})} + \|\frac{\partial u}{\partial x}\|_{L^\infty(\Omega\times \mathcal{T})})$ and $C_1$ is defined in Proposition \ref{ProAsympErr}. $M_0$ represents the set of  bounded functions, i.e.
\begin{equation*}
M_0 := \left\{ f: f(x) \in[C_l, C_u], \forall x\in\Omega  \right\}.
\end{equation*}

It can be clearly shown that both $f^\delta$ and $f^*$ belong to $Q^{1,2}_{ad}$. Hence, once the approximate solution $f^\delta$ is obtained, we can calculate the a posteriori error of $f^\delta$ according to the optimization problem
\begin{equation}
\label{posterioriErr}
\max_{f\in Q^{1,2}_{ad}} \|f-f^\delta\|^2_{L^p(\Omega)},
\end{equation}
which is well-posed according to the Bolzano–Weierstrass theorem – note that $Q^{1,2}_{ad}$ are compact sets in $L^2(\Omega)$. The $n$-dimensional analog of \eqref{posterioriErr} is
\begin{equation}
\label{posterioriErr2}
\bar{\Delta}_1 = \max_{\mathbf{f}\in \mathbf{Q}^{1,2}_{ad}} \|\mathbf{f}-\mathbf{f}^\delta\|^2, \quad \Delta_1 = \sqrt{\bar{\Delta}_1}/\|\mathbf{f}^\delta\|,
\end{equation}
where we call the number $\Delta_1$ the relative a posteriori error for the reconstructed source function $\mathbf{f}^\delta$. Here, $\mathbf{f}, \mathbf{f}^\delta \in \mathbb{R}^n$, with $\mathbf{f}^\delta$ being the $n$-dimensional projection of $f^\delta$, and the $n$-dimensional analogs of the admissible sets $Q^{1,2}_{ad}$ are defined as follows.

The constraint for bounded monotonic functions:
\begin{equation}
\label{Q1}
\mathbf{Q}^{1}_{ad} = \{\mathbf{f} \in \mathbb{R}^n: \mathbf{A}_1 \mathbf{f} \leq \mathbf{b}_1 \},
\end{equation}
where
$$
\mathbf{A_1} = \left(\begin{array}{cccccc}
1 & -1 & 0 & \ldots & 0 & 0 \\
0 & 1 & -1 & \ldots & 0 & 0 \\
\vdots & \vdots & \vdots & \ddots & \vdots & \vdots  \\
0 & 0 & 0 & \ldots & 1 & -1 \\
\hline 1 & 0 & 0 & \ldots & 0  & 0  \\
\vdots & \vdots & \vdots & \ddots & \vdots & \vdots \\
0 & 0 & 0 & \ldots & 0  & 1\\
\hline -1 & 0 &  0 & \ldots & 0 & 0 \\
\vdots & \vdots &  \vdots & \ddots & \vdots  & \vdots \\
0 & 0 & 0 &  \ldots & 0 & -1
\end{array}\right),  \quad
\mathbf{b}_1= \left(\begin{array}{c}
0 \\
0 \\
\vdots \\
0 \\
\hline C_u \\
\vdots \\
C_u \\
\hline  C_l \\

\vdots \\
C_l
\end{array}\right).
$$

The constraint for bounded convex functions:
\begin{equation}
\label{Q2}
\mathbf{Q}^{2}_{ad} = \{\mathbf{f} \in \mathbb{R}^n: \mathbf{A}_2 \mathbf{f} \leq \mathbf{b}_2 \},
\end{equation}
where (the grid $\Theta $ is assumed to be uniform)
$$
\mathbf{A_2} = \left(\begin{array}{cccccc}
1 & -2 & 1 & \ldots & 0 & 0 \\
0 & 1 & -2 & \ldots & 0 & 0 \\
\vdots & \vdots & \vdots & \ddots & \vdots & \vdots  \\
0 & 0 & 0 & \ldots & -2 & 1 \\
\hline 1 & 0 & 0 & \ldots & 0  & 0  \\
\vdots & \vdots & \vdots & \ddots & \vdots & \vdots \\
0 & 0 & 0 & \ldots & 0  & 1\\
\hline -1 & 0 &  0 & \ldots & 0 & 0 \\
\vdots & \vdots &  \vdots & \ddots & \vdots  & \vdots \\
0 & 0 & 0 &  \ldots & 0 & -1
\end{array}\right),  \quad
\mathbf{b}_2= \left(\begin{array}{c}
0 \\
0 \\
\vdots \\
0 \\
\hline C_u \\
\vdots \\
C_u \\
\hline  C_l \\

\vdots \\
C_l
\end{array}\right).
$$

To obtain the pointwise error estimate for the reconstructed source function $f^\delta$, we construct the upper solution $f^{up}$ and the lower solution $f^{low}$ so that, for both $f^\delta$ and $f^*$, the following holds for all $x\in \Omega$:
\begin{equation} \label{UpandLowSol}
f^{low}(x) \leq  f^\delta(x),\  f^*(x)  \leq f^{up}(x).
\end{equation}

Once the upper and lower solutions are constructed, the pointwise error estimate can be calculated through
\begin{equation}
\label{pointwiseErr}
\Delta_2(x) = f^{up}(x)- f^{low}(x).
\end{equation}

To do this, for $i=0, \cdots, n$, let
\begin{equation}
\label{minmax_points}
\mathbf{f}_{i}^{low}=\inf \left\{\mathbf{f}_{i}: \mathbf{f} \in \mathbf{Q}^{1,2}_{ad} \right\}, \ \mathbf{f}_{i}^{up}=\sup \left\{\mathbf{f}_{i}: \mathbf{f} \in \mathbf{Q}^{1,2}_{ad} \right\}.
\end{equation}

We construct the upper and lower solutions with points from \eqref{minmax_points} (a similar idea can be found in \cite{TitarenkoYagola2002}).

For the monotonic functions, it is clear that $\mathbf{f}_{i}^{low} \leq \mathbf{f}_{i+1}^{low}, \mathbf{f}_{i}^{up} \leq \mathbf{f}_{i+1}^{low}, \mathbf{f}_{i}^{low} \leq \mathbf{f}_{i}^{up}, i=0, \cdots ,n-1$. Hence, the lower and upper solutions can be constructed as follows (see Fig. \ref{fig:monotonicUPLOW} for the visualization):
\begin{align}
f^{low}(x) &=\left\{\begin{array}{ll}
\mathbf{f}_{0}^{low}, & x \in\left[x_{0}, x_{1}\right], \label{monotonicLOW} \\
\mathbf{f}_{i}^{low}, & x \in\left(x_{i}, x_{i+1}\right],
\end{array} \quad i=1, \cdots, n-1, \right. \\
f^{up}(x) &= \left\{\begin{array}{ll}
\mathbf{f}_{i+1}^{up}, & x \in\left[x_{i}, x_{i+1}\right), \label{monotonicUP} \\
\mathbf{f}_{n}^{up}, & x \in\left[x_{n-1}, x_{n}\right],
\end{array} \quad i=0, \cdots, n-2. \right.
\end{align}

For the set of convex functions $M_2$, the upper and lower solutions are constructed according to the following equations (see Fig. \ref{fig:convexUPLOW} for the visualization):
\begin{align}\label{convexLOW}
f^{low}(x)=\frac{\mathbf{f}_{i+1}^{low}-\mathbf{f}_{i}^{low}}{x_{i+1}-x_{i}} x+\frac{\mathbf{f}_{i}^{low} x_{i+1}-\mathbf{f}_{i+1}^{low} x_{i}}{x_{i+1}-x_{i}}, \quad x \in\left[x_{i}, x_{i+1}\right], \quad i=0, \cdots, n-1,
\end{align}
\begin{align}\label{convexUP}
f^{up}(x)= \begin{cases} \displaystyle  \frac{\mathbf{f}^{up}_{1}-\mathbf{f}^{up}_{0}}{x_{1}-x_{0}} x+\frac{\mathbf{f}^{up}_{0} x_{1}-\mathbf{f}^{up}_{1} x_{0}}{x_{1}-x_{0}}, & x \in\left[x_{0}, x_{1}\right], \\
\displaystyle \frac{\mathbf{f}^{up}_{i}-\mathbf{f}^{low}_{i-1}}{x_{i}-x_{i-1}} x+\frac{\mathbf{f}^{low}_{i-1} x_{i}-\mathbf{f}^{up}_{i} x_{i-1}}{x_{i}-x_{i-1}}, & x \in\left[x_{i}, x_{i}^{\pm}\right], \\
\displaystyle \frac{\mathbf{f}^{low}_{i+2}-\mathbf{f}^{up}_{i+1}}{x_{i+2}-x_{i+1}} x+\frac{\mathbf{f}^{up}_{i+1} x_{i+2}-\mathbf{f}^{low}_{i+2} x_{i+1}}{x_{i+2}-x_{i+1}}, & x \in\left[x_{i}^{\pm}, x_{i+1}\right], \\
\displaystyle \frac{\mathbf{f}^{up}_{i}-\mathbf{f}^{up}_{i-1}}{x_{i}-x_{i-1}} x+\frac{\mathbf{f}^{up}_{i-1} x_{i}-\mathbf{f}^{up}_{i} x_{i-1}}{x_{i}-x_{i-1}}, & x \in\left[x_{n-1}, x_{n}\right], \end{cases}
\end{align}
where $i=1, \cdots, n-2$ and
$$
x_{i}^{\pm}=\frac{\left(\mathbf{f}^{low}_{i-1}x_{i}-\mathbf{f}^{up}_{i} x_{i-1}\right)\left(x_{i+2}-x_{i+1}\right)-\left(\mathbf{f}^{up}_{i+1} x_{i+2}-\mathbf{f}^{low}_{i+2} x_{i+1}\right)\left(x_{i+2}-x_{i+1}\right)}{\left(\mathbf{f}^{low}_{i-1}-\mathbf{f}^{up}_{i}\right)\left(x_{i+2}-x_{i+1}\right)-\left(\mathbf{f}^{up}_{i+1}-\mathbf{f}^{low}_{i-1}\right)\left(x_{i}-x_{i-1}\right)}.
$$

Based on the above analysis, we construct an efficient regularization algorithm for the nonlinear inverse source problem (\textbf{IP}), as shown below.

\begin{algorithm}[!htb]
\caption{Asymptotic expansion-regularization (AER) algorithm for (\textbf{IP}).}
\label{alg:Framwork}
\begin{algorithmic}[1]
  \item[] \begin{itemize}  \item[\%]  Calculation of the approximate source function. \end{itemize}
\If{The full measurement data $\{u^\delta_i, w^\delta_i\}^n_{i=0}$ are given}
    \State break;
\Else
    \If{ Only the measurements $\{u^\delta_i\}^n_{i=0}$ are provided}
        \State Construct the smoothed data $\{u^\varepsilon, \frac{\partial u^\varepsilon}{\partial x}\}$ by solving \eqref{uAlpha};
    \EndIf
\EndIf
\State Calculate the approximate source function $f^\delta$ using formula \eqref{uAlpha1}.
\begin{itemize}
\item[\%]  Error estimation.
\end{itemize}
\State Calculate the a posteriori error $\Delta_1$ of the obtained approximate source function by solving problem  \eqref{posterioriErr2}.

\State Find the lower $f^{low}$ and upper $f^{up}$ solutions using formulas  \eqref{monotonicLOW}–\eqref{convexUP}.

\State Calculate the pointwise error estimate $\Delta_2(x)$ using formula \eqref{pointwiseErr}.
\end{algorithmic}
\end{algorithm}


\section{Derivation and proofs of main results} \label{derivationAndProofs}

We follow the idea in \cite{b6} and consider a solution in the form of a moving front, which at each moment of time $t$ is localized in a neighborhood of some point $x_{t.p.}(t,\mu)\in \Omega$ to the left and to the right of it; a narrow moving transition layer is observed in the indicated vicinity.

\subsection{Construction of the asymptotic solution}\label{constructionOfAsymptotic}

The asymptotic solution of problem \eqref{equat1} will be constructed in the following form:
\begin{align}
\label{u}
U=\begin{cases}
U^{l}, &  (x,t)\in  \bar{\Omega}^{l} \times \bar{\mathcal{T}}:= \{(x,t) \in \mathbb{R}^2: x\in [0,x_{t.p.}(t,\mu)],  t \in \bar{\mathcal{T}} \}  ,\\
U^{r}, &    (x,t)\in  \bar{\Omega}^{r} \times \bar{\mathcal{T}} := \{(x,t) \in \mathbb{R}^2: x\in [x_{t.p.}(t,\mu),1], t \in \bar{\mathcal{T}} \} ,
\end{cases}
\end{align}
where $\bar{\Omega}^{l} \times \bar{\mathcal{T}}$ and  $\bar{\Omega}^{r} \times \bar{\mathcal{T}}$ -- are the regions to the left and right, respectively, of the point  $x_{t.p.}(t,\mu)$, and the functions  $U^{l}$ and $U^{r}$ have the following form:
\begin{align}
\label{ulr}
U^{l}=\bar{u}^{l}(x,\mu)+Q^{l}(\xi,t,\mu),\qquad
U^{r}=\bar{u}^{r}(x,\mu)+Q^{r}(\xi,t,\mu),
\end{align}
where $\bar{u}^{l,r}(x,\mu)$ -- are regular functions describing the solution away from the point $x_{t.p.}(t,\mu)$, and the functions $Q^{l,r}(\xi,t,\mu)$ describe the transition layer near the point $x_{t.p.}(t,\mu)$, with the variable $\xi$ defined as
\begin{align}
\label{xi}
\xi := \frac{x-x_{t.p.}(t,\mu)}{\mu} \begin{cases}
\quad \leq 0, &  \bar{\Omega}^{l} \times \bar{\mathcal{T}},\\
\quad \geq 0, &  \bar{\Omega}^{r} \times \bar{\mathcal{T}}.
\end{cases}
\end{align}

Note that, from Assumption \ref{A4}, the function $u(x)$ has a transition layer between the levels $\varphi^{l}(x)$ and $\varphi^{r}(x)$  in the vicinity of the point $x_{0}^{*} \in \Omega$.
We take $x_{t.p.}(0,\mu)=x_{0}^{*}$ and look for the coordinate and the velocity of the transitional layer in the form \footnote{Hereinafter, functions with subscript $0$ will be called the zero approximation of a function expandable in powers of $\mu$. In this paper, we mainly focus on the zero approximation for the asymptotic solution.}
\begin{equation} \label{expantionxtp}
 x_{t.p.}(t,\mu)=x_{0}(t)+\mu \cdot x_{1}(t)+\cdots,
  \end{equation}
\begin{equation}\label{expantionv}
 v_{t.p.}(t,\mu)=v_{0}(t)+\mu \cdot v_{1}(t)+\cdots,
  \end{equation}
where $v_{i}(t)=\displaystyle dx_{i} / dt$.

 The functions $\bar{u}^{l,r}(x,\mu)$ and $Q^{l,r}(\xi,t,\mu)$ in \eqref{ulr} are constructed as expansions in powers of  $\mu$:
\begin{equation} \label{regularexpansions}
\bar{u}^{l,r}(x,\mu)=\bar{u}_{0}^{l,r}(x)+\mu \bar{u}_{1}^{l,r}(x)+\cdots+\mu^{n} \bar{u}_{n}^{l,r}(x)+\cdots,
\end{equation}
\begin{equation} \label{transitionexpansions}
 Q^{l,r}(\xi,t,\mu)=Q_{0}^{l,r}(\xi,t)+\mu Q_{1}^{l,r}(\xi,t)+\cdots +\mu^{n}Q_{n}^{l,r}(\xi,t)+\cdots .
\end{equation}

For the regular asymptotic part we substitute expansions \eqref{regularexpansions} into the stationary equations, i.e. $du/dt =0$:
\begin{equation} \label{regularpartfirstorder}
\displaystyle \mu \frac{\partial^{2} \bar{u}^{l,r}}{\partial x^{2}}=-k\bar{u}^{l,r} \frac{\partial \bar{u}^{l,r}}{\partial x}+ f(x).
\end{equation}

Equating coefficients at $\mu =0$ in \eqref{regularpartfirstorder}, we obtain the degenerate equations for the functions $\varphi^{l} (x)$ and $\varphi^{r} (x)$ with the left and right boundary conditions, respectively:
\begin{align}
\label{LeftRightSolution}
 \begin{cases}
 \displaystyle -k \varphi^{l} (x) \frac{d\varphi^{l} (x)}{dx}+f(x)=0, ~ x\in \bar{\Omega}, \\
\displaystyle \varphi^{l}(0) =u^{l},
\end{cases}
\begin{cases}
\displaystyle -k \varphi^{r} (x) \frac{d\varphi^{r} (x)}{dx}+f(x)=0, ~ x\in \bar{\Omega}, \\
\displaystyle \varphi^{r}(1)=u^{r},
\end{cases}
\end{align}
where the main terms of the regular part of $U$, i.e. \eqref{u}, are defined as
\begin{align}
\label{u0regu}
\bar{u}_{0}(x)= \begin{cases}
\bar{u}_{0}^{l}(x)=\varphi^{l}(x), \quad x\in [0,x_{t.p.}(t,\mu)],\\
\bar{u}_{0}^{r}(x)=\varphi^{r}(x), \quad x\in [x_{t.p.}(t,\mu),1].
\end{cases}
\end{align}

Taking into account the fact that
\begin{equation} \label{operatorequat}
\mu \frac{\partial^2}{\partial x^2}- \frac{\partial}{\partial t}=\frac{1}{\mu}\frac{\partial^2}{\partial \xi^2}+\frac{1}{\mu} v_{t.p}(t,\mu)\frac{\partial}{\partial \xi} - \frac{\partial}{\partial t}, \quad \frac{\partial}{\partial x}= \frac{1}{\mu} \frac{\partial}{\partial \xi},
\end{equation}
we substitute the expansions \eqref{ulr} into \eqref{equat1}. Then, subtracting the regular part from the result, we obtain the equation for transition-layer functions:
\begin{multline} \label{TLmaineq}
 \frac{1}{\mu}\frac{\partial^{2}Q^{l,r}}{\partial\xi^{2}}+\frac{1}{\mu}v_{t.p.}(t, \mu)\frac{\partial Q^{l,r}}{\partial\xi}-\frac{\partial Q^{l,r}}{\partial t}
 =- \frac{k}{\mu} \Big( \left(\bar{u}^{l,r}(\mu\xi+x_{t.p}(t, \mu), \mu)+Q^{l,r}(\xi, t, \mu)\right) (\frac{\partial Q^{l,r}}{\partial\xi}+\frac{ \partial \bar{u}^{l,r}(\mu\xi+x_{t.p}(t, \mu), \mu)}{\partial \xi})\\
 - \left(\bar{u}^{l,r}(\mu\xi+x_{t.p}(t, \mu), \mu)\right)\frac{\partial \bar{u}^{l,r}(\mu\xi+x_{t.p}(t, \mu), \mu)}{\partial \xi} \Big).
\end{multline}
Substituting the expansions \eqref{expantionxtp}–\eqref{transitionexpansions} into \eqref{TLmaineq}  and equating the coefficients at $\mu^{-1}$, we obtain
\begin{align} \label{equat21}
 \begin{cases}
 \displaystyle \frac{\partial^{2}Q_{0}^{l,r}}{\partial\xi^{2}}+ \bigg(k \Big(\varphi^{l,r} (x_{0}(t))+Q_{0}^{l,r} \Big)+ v_{0} (t) \bigg) \frac{\partial Q_{0}^{l,r}}{\partial\xi}=0,\\
Q_{0}^{l}(0,t)+\varphi^{l}(x_{0}(t))=\varphi(x_{0}(t)), \\
Q_{0}^{r}(0,t)+\varphi^{r}(x_{0}(t))=\varphi(x_{0}(t)), \\
Q_{0}^{l}(\xi,t)\rightarrow 0 \ \text{for} \ \xi\rightarrow-\infty,\\
Q_{0}^{r}(\xi,t)\rightarrow 0 \ \text{for} \ \xi\rightarrow+\infty,
\end{cases}
\end{align}
where
\begin{equation} \label{varphi}
\displaystyle \varphi(x_{0}(t)) :=\frac{1}{2} \left( \varphi^{l}(x_{0}(t))+\varphi^{r}(x_{0}(t)) \right).
\end{equation}

To study the zero approximation of $x_{t.p.}(t,\mu)$, i.e. $x_0(t)$, we introduce the auxiliary function
\begin{align} \label{auxilaryfunction}
 \tilde{u}=
\begin{cases}
\varphi^{l}(x_{0}(t))+Q_{0}^{l}(\xi,t), \quad & \xi \leq 0, \quad t \in \bar{\mathcal{T}}, \\
\varphi(x_{0}(t)), \quad & \xi = 0, \quad t \in \bar{\mathcal{T}}, \\
 \varphi^{r}(x_{0}(t))+Q_{0}^{r}(\xi,t), \quad &\xi \geq 0, \quad t \in \bar{\mathcal{T}}. \\
\end{cases}
\end{align}
It is clear that $\tilde{u} \in (\varphi^{l}(x_{0}(t),\varphi^{r}(x_{0}(t)))$. We rewrite problem \eqref{equat21} in the following form:
\begin{align} \label{replacement}
 \begin{cases}
\displaystyle \frac{\partial^2 \tilde{u} }{\partial\xi^{2}} = (-k \tilde{u}-v_{0}(t)) \frac{\partial \tilde{u} }{\partial\xi}, \\
\displaystyle  \tilde{u} (0,t)= \varphi(x_{0}(t)), ~ \tilde{u} (-\infty ,t) = \varphi^{l}(x_{0}(t)), ~ \tilde{u} (+ \infty ,t) = \varphi^{r}(x_{0}(t)).
 \end{cases}
 \end{align}

Let $\displaystyle \frac{\partial \tilde{u} }{\partial\xi} = g(\tilde{u}), \frac{\partial^2 \tilde{u} }{\partial\xi^{2}} = \frac{\partial g(\tilde{u}) }{\partial \tilde{u} } g(\tilde{u}) $; then, \eqref{replacement} is transformed into $\displaystyle \frac{\partial g(\tilde{u}) }{\partial \tilde{u} } = -k \tilde{u}-v_{0}(t)$, from which we can deduce that
\begin{align} \label{replacement2}
\frac{\partial \tilde{u} }{\partial\xi} =
 \begin{cases}
\displaystyle \int_{\varphi^{l}(x_{0}(t))}^{\tilde{u}} (-k u - v_{0} (t)) du, \quad  \xi \leq 0, \\ \\
\displaystyle  \int_{\varphi^{r}(x_{0}(t))}^{\tilde{u}}( -k u - v_{0} (t)) du, \quad  \xi \geq 0.
 \end{cases}
 \end{align}

From the zeroth-order $C^1$-matching condition, i.e.
\begin{equation}
\frac{\partial \tilde{u}}{\partial \xi } \Big \vert_{\begin{subarray}{l} \xi=-0 \end{subarray}}=  \frac{\partial \tilde{u}}{\partial \xi } \Big \vert_{\begin{subarray}{l} \xi=+0 \end{subarray}}, \quad t \in \bar{\mathcal{T}}, \label{matching}
\end{equation}
we obtain, for the zero approximation,
\begin{equation}\label{EqV0}
\int_{\varphi^{l}(x_{0}(t))}^{\varphi^{r}(x_{0}(t))} (-ku-v_0 (t)) du =0,
\end{equation}
where $v_0 (t) = dx_0(t) / dt$. Solving \eqref{EqV0}, we obtain
\begin{equation} \label{speed}
v_{0}(t) = -\frac{k}{2} (\varphi^{r}(x_{0}(t))+\varphi^{l}(x_{0}(t))).
\end{equation}
By using the explicit formula for $\varphi^{l,r}(x_{0}(t))$ (cf. \eqref{eq7}), we obtain the equation that determines the location of the transition layer in the zero approximation $x_{0}(t) $ for every $t$:
\begin{align}
\begin{cases} \label{eq10}
\displaystyle \frac{dx_0(t)}{dt}=\frac{k}{2} \left( \sqrt{\frac{2}{k} \int_0^{x_{0}(t)} f(s )  d s + (u^{l})^2} -\sqrt{  (u^{r})^2 - \frac{2}{k} \int_{x_{0}(t)}^1  f(s ) d s} \right),\\
x_0(0)=x_{0}^{*}\in \Omega.
\end{cases}
\end{align}

\begin{remark}
\label{remarkx0}
Empirically, we found that the quantity $ \vert \varphi^{l} (x) \vert - \vert \varphi^{r} (x) \vert$ keeps the sign for all $x \in \Omega$. In the case when $ \vert \varphi^{l} (x) \vert >\vert \varphi^{r} (x) \vert $,  the right-hand side of the ordinary differential equation \eqref{eq10} is always positive, which means that the solution of \eqref{eq10} is increasing, and hence Assumption \ref{A3} is reduced to check the condition $x_0(T)<1$ for only one time point $t=T$. Alternatively, if, for any $x \in \Omega $, $ \vert \varphi^{l} (x) \vert < \vert \varphi^{r} (x) \vert $ holds, the solution of \eqref{eq10} is decreasing and Assumption \ref{A3} can be simplified as the positivity of $x_0(T)$.
\end{remark}

Integrating the right-hand side of \eqref{replacement2}, we obtain
\begin{align} \label{derivativetildeu}
 \frac{\partial \tilde{u} }{\partial\xi} =
\begin{cases}
\displaystyle   \frac{k}{2}( \varphi^{l}(x_{0}(t))^2-\tilde{u}^2)  - \left( \tilde{u} -\varphi^{l}(x_{0}(t))\right) v_0(t), &  \xi \leq 0, \\
\displaystyle  \frac{k}{2}( \varphi^{r}(x_{0}(t))^2-\tilde{u}^2)  - \left( \tilde{u} -\varphi^{r}(x_{0}(t))\right)v_0(t), &  \xi \geq 0.
 \end{cases}
 \end{align}

From equations \eqref{speed} and \eqref{derivativetildeu} and the definition of $\tilde{u}$ in \eqref{auxilaryfunction}, we can write the functions $Q_{0}^{l,r}(\xi, t)$ in the explicit form, in which $ x_0 (t) $ is a parameter:
\begin{equation} \label{zeroordertransitionfunc}
\displaystyle Q_{0}^{l,r}( \xi, t)= \frac{ -2P^{l,r}(x_0(t))}{  \exp \left(\xi k P^{l,r}(x_0(t)) \right)+ 1} ,
\end{equation}
where
\begin{align*}
& P^{l}(x_0(t))=\frac{1}{2}\left(\varphi^{l}(x_{0}(t))-\varphi^{r}(x_{0}(t)) \right), \\
& P^{r}(x_0(t))=\frac{1}{2}\left(\varphi^{r}(x_{0}(t))-\varphi^{l}(x_{0}(t))\right)=-P^{l}(x_0(t)).
\end{align*}

According to the boundary conditions at points $\xi=\pm\infty$ in \eqref{equat21}, we conclude that $P^{l}(x_0(t))<0$ and $P^{r}(x_0(t))>0$ for all $t>0$.

Consequently, the functions $Q_{0}^{l,r}(\xi, t)$ are exponentially decreasing with $\xi \rightarrow \mp \infty $ and have the exponential estimates (see, e.g., \cite{Vasileva1998ContrastSI,Butuzov1997ASYMPTOTICTO})
\begin{equation}\label{equat22}
\underline{C} e^{\underline{\kappa}\xi} \leq \lvert Q_{0}^{l}(\xi, t) \rvert \leq \bar{C}e^{\bar{\kappa}\xi} , \quad \xi\leq 0, \quad t\in \bar{\mathcal{T}},
 \end{equation}
\begin{equation}\label{equat23}
\underline{C} e^{-\underline{\kappa}\xi} \leq \lvert Q_{0}^{r}(\xi, t) \rvert \leq \bar{C}e^{-\bar{\kappa}\xi}, \quad \xi\geq 0, \quad t\in \bar{\mathcal{T}},
 \end{equation}
where $\underline{C}, \bar{C}$ and $\underline{\kappa}, \bar{\kappa}$ -- are four positive constants, and, more precisely, $$\underline{C} := \frac{1}{2}  \inf_{t \in [0,T]} \left( \varphi^{r}(x_{0}(t))-  \varphi^{l}(x_{0}(t)) \right), \quad \bar{C} := \frac{1}{2}  \sup_{t \in [0,T]} \left( \varphi^{r}(x_{0}(t))-  \varphi^{l}(x_{0}(t)) \right).$$

From the boundary conditions of \eqref{equat21} and Assumption \ref{A1}, we deduce that
\begin{equation*}
\lvert Q_{0}^{l,r}(0,t)\rvert=\lvert \varphi(x_{0}(t))-\varphi^{l,r}(x_{0}(t)) \rvert = \frac{1}{2} (\varphi^{r}(x_{0}(t)) - \varphi^{l}(x_{0}(t))) \geq \frac{1}{2} \left( u^r- u^l \right) > \mu^2,
\end{equation*}
and $\lvert Q_{0}^{l,r}(\xi,t) \rvert \rightarrow 0$  for $\xi \rightarrow \mp \infty$. Since $\mu >0 $ is a fixed number, $\lvert Q_{0}^{l,r}(\xi,t) \rvert$ are decreasing functions and  $\xi=(x-x_{t.p.}(t,\mu))/ \mu $, then there exist  $\hat{x}^{l,r}(t,\mu)$ for which on the intervals $ [0, \hat{x}^{l}(t,\mu)]$ and $[ \hat{x}^{r}(t,\mu),1]$ we have  $\lvert Q_{0}^{l,r}(\xi, t) \rvert \leq \mu^2 $ respectively for every $t$; and at the points $\hat{x}^{l,r}(t,\mu)$:
\begin{equation}\label{Qlr-mu}
\lvert Q_{0}^{l}(\xi( \hat{x}^{l}(t,\mu)), t) \rvert = \mu^2, \quad \lvert Q_{0}^{r}(\xi( \hat{x}^{r}(t,\mu)), t) \rvert = \mu^2.
 \end{equation}

Now, we define the width of transition layer $\Delta x(t,\mu)$ and the point at the middle of the transition layer, $x_{t.p.}(t,\mu)$, with
\begin{equation}
\label{transitionVariables}
\Delta x(t,\mu) = \hat{x}^{r}(t,\mu) -\hat{x}^{l}(t,\mu), \quad x_{t.p.}(t,\mu) = (\hat{x}^{r}(t,\mu) + \hat{x}^{l}(t,\mu))/2.
\end{equation}

From equations \eqref{equat22}–\eqref{equat23} and \eqref{Qlr-mu} for  $x =\hat{x}^{l,r}(t)$  we obtain
\begin{equation}
\underline{C} e^{-\frac{\underline{\kappa}}{2\mu} \Delta x} \leq \mu^2 \leq \bar{C}e^{-\frac{\bar{\kappa}}{2\mu} \Delta x},
\end{equation}
from which $ \Delta x$ can be estimated as
\begin{equation} \label{DeltaX}
\frac{2\mu}{\underline{\kappa}} \ln \frac{\underline{C}}{\mu^2} \leq \Delta x \leq \frac{2\mu}{\bar{\kappa}} \ln \frac{\bar{C}}{\mu^2}, \text{~i.e.~} \Delta x \sim \mu \lvert \ln \mu \rvert .
\end{equation}

We also write the first-order asymptotic approximation functions. Equating the coefficients at $\mu^1$ in \eqref{regularpartfirstorder}, we obtain the following equations:
\begin{align}
\begin{split}
&\displaystyle  \varphi^{l,r}(x)\frac{d \bar{u}_{1}^{l,r}}{d x}+\bar{u}_{1}^{l,r}\frac{d \varphi^{l,r}(x) }{d x}=- \frac{1}{k} \frac{d^2 \varphi^{l,r}(x)}{d x^2}, \\
&\bar{u}_{1}^{l}(0,t)=0,\quad \bar{u}_{1}^{r}(1,t)=0.
\end{split}
\end{align}

The solutions to these problems can be written explicitly:
\begin{align}
\begin{split}
\bar{u}_{1}^{l}(x,t) = \exp \left( \int_{0}^{x} -W(s) ds \right) \int_{0}^{x}- \exp \left( \int_{0}^{s'} W(s) ds \right) Y(s') ds', \\
\bar{u}_{1}^{r}(x,t) = \exp \left( \int_{x}^{1} W(s) ds \right) \int_{x}^{1} \exp \left( \int_{s'}^{1} -W(s) ds \right) Y(s') ds',
\end{split}
\end{align}
where $W(x)=\displaystyle \frac{1}{\varphi^{l,r}(x)}\frac{d \varphi^{l,r}(x) }{d x}$ and $Y(x)= \displaystyle \frac{1}{k\varphi^{l,r}(x)}\frac{d^2 \varphi^{l,r}(x) }{d x^2} $.

After substituting expansions
\begin{align*}
\displaystyle \bar{u}^{l,r}(\mu\xi+x_{0}(t)+\mu x_1 (t),\mu) =  u_0(x_0(t))+\mu \left((\xi+x_1) \frac{d u_0 }{d x}(x_0(t))  + u_1(x_0(t)) \right)+\mathcal{O}(\mu^2),
\end{align*}
and expansions \eqref{expantionxtp},\eqref{expantionv}, and \eqref{transitionexpansions} into \eqref{TLmaineq}  and equating the coefficients at $\mu^0$, we obtain equations for the first-order transition-layer functions:
\begin{align*}
\begin{split}
 &\frac{\partial^{2} Q_{1}^{l,r} }{\partial \xi^2}+ \left(k(\varphi^{l,r}(x_0(t))+Q_{0}^{l,r} ) +v_0(t) \right)\frac{\partial Q_{1}^{l,r} }{\partial \xi}  +  k  Q_{1}^{l,r}  \Upsilon^{l,r}(\xi,t) \\
 &=   \left(  -k(x_1(t) \frac{d \varphi^{l,r} }{d x}(x_0(t))  + \bar{u}_{1}^{l,r}(x_0(t))) - v_1 (t) \right) \Upsilon^{l,r} + r_{1}^{l,r}(\xi,t) :=H_1^{l,r} (\xi,t),
\end{split}
\end{align*}
where
\begin{equation*} 
\Upsilon^l(\xi,t) = \frac{\partial \tilde{u}}{\partial \xi} (\xi,t), \ \xi \leq 0 , \quad \Upsilon^r(\xi,t) = \frac{\partial \tilde{u}}{\partial \xi} (\xi,t), \ \xi \geq 0,
\end{equation*}
and
\begin{align*}
\begin{split}
r_{1}^{l,r}(\xi,t) = -k \left( \left(   \xi  \frac{d \varphi^{l,r} }{d x}(x_0(t))  \right)\Upsilon^{l,r}+Q_{0}^{l,r} \frac{d   \varphi^{l,r} }{ d x}(x_0(t)) \right)+\frac{\partial Q_{0}^{l,r}}{\partial t}.
\end{split}
\end{align*}

Taking into account the initial conditions in \eqref{equat21}, we derive additional conditions for the functions $Q_{1}^{l,r} (\xi,t)$:
\begin{align*}
Q_{1}^{l,r} (0,t)=-\bar{u}_{1}^{l,r}(x_0(t))- x_1(t) \frac{d \varphi^{l,r} }{d x}(x_0(t))\equiv p_{1}^{l,r} (t), \\
Q_{1}^{l} ( \xi ,t) \rightarrow 0 \ \text{for} \  \xi \rightarrow -\infty, \quad Q_{1}^{r} (\xi,t) \rightarrow 0 \ \text{for} \ \xi \rightarrow +\infty.
\end{align*}

Functions $Q_{1}^{l,r}(\xi,t)$ can also be written explicitly:
\begin{equation} \label{Q1function}
Q_{1}^{l,r} (\xi,t)=z^{l,r}(\xi,t) \left( p_{1}^{l,r} (t)   - \int_{0}^{\xi} \frac{1}{z^{l,r}(s,t)} \int_{s}^{\mp \infty} H_{1}^{l,r} (\eta,t) d\eta ds \right),
\end{equation}
where $ \displaystyle z^{l,r}(\xi,t)= \left( \Upsilon^{l,r}(0,t) \right)^{-1} \Upsilon^{l,r}(\xi,t)$.  From \eqref{Q1function} we deduce that
\begin{align*}
\frac{\partial Q_{1}^{l,r}}{\partial \xi} (0,t) =\left(-k\varphi^{l,r}(x_0(t)) -v_0(t) \right) p_{1}^{l,r} (t)
+\frac{dx_1(t)}{dt} \left(\varphi^{l,r}(x_0(t))-\varphi (x_0(t))\right)- \int_{0}^{\mp \infty} r_{1}^{l,r} (\eta,t) d\eta.
\end{align*}

It can be clearly shown that the functions $Q_{1}^{l,r} (\xi,t) $ satisfy exponential estimates of the type \eqref{equat22}, \eqref{equat23}. From the first-order $C^1$-matching condition
\begin{equation} \label{matchingfirstord}
\frac{\partial Q_{1}^{l}}{\partial \xi }(0,t) + \frac{d \varphi^{l}}{d x} (x_0 (t))=\frac{\partial Q_{1}^{r}}{\partial \xi }(0,t) + \frac{d \varphi^{r}}{d x} (x_0 (t)), \quad t \in \bar{\mathcal{T}},
\end{equation}
we obtain the equation that determines $x_1(t)$:

\begin{equation} \label{eqforx1}
(\varphi^l (x_0(t))-\varphi^r (x_0(t)))\frac{d x_1 (t)}{d t}+\Phi_1(t) x_1(t)=\Phi_2(t),
\end{equation}
where
\begin{align*}
 \Phi_1(t)=&  \frac{k}{2}\left(\varphi^{r}(x_0(t))-\varphi^{l}(x_0(t))  \right) \left( \frac{d \varphi^{r} }{d x}(x_0(t))+ \frac{d \varphi^{l} }{d x}(x_0(t)) \right), \\  \Phi_2(t)=&  \frac{k}{2} \left(  \bar{u}_{1}^{l}(x_0(t))   + \bar{u}_{1}^{r}(x_0(t))  \right) \left(\varphi^{l}(x_0(t)) - \varphi^{r}(x_0(t)) \right)  -\frac{d \varphi^{l} }{d x} (x_0(t)) +\frac{d \varphi^{r} }{d x} (x_0(t))+ \int_{0}^{- \infty} r_{1}^{l} (\eta,t) d\eta -\int_{0}^{+ \infty} r_{1}^{r} (\eta,t) d\eta.
\end{align*}

Given the initial condition in \eqref{eq10} and that $x_{t.p.}(0,\mu)=x_{0}^{*}$, we solve \eqref{eqforx1} with initial condition $x_1(0)=0$, so we can find $x_1(t)$  in explicit form:
\begin{align} \label{x1explicit}
x_1(t)= \mathrm{exp}  \left( - \int_{0}^{t} \frac{\Phi_1 (s)}{\varphi^l (x_0(t))-\varphi^r (x_0(t))} ds     \right)  \int_{0}^{t}   \frac{ \mathrm{exp} \left(  \int_{0}^{\eta} \frac{\Phi_1 (s)}{\varphi^l (x_0(t))-\varphi^r (x_0(t))} ds     \right) \Phi_2 (\eta) }{\varphi^l (x_0(t))-\varphi^r (x_0(t))} d \eta.
\end{align}

In a similar way to \eqref{regularpartfirstorder}–\eqref{x1explicit}, we can obtain approximation terms of solution $u$ up to $n$ order, i.e. formula \eqref{asymptoticnorder}.

Moreover, the approximation terms of $x_{t.p}(t,\mu)$ in \eqref{expantionxtp} up to order $n$ can be written as
\begin{equation} \label{XtpOfOrderN}
X_{n}(t,\mu)=\sum_{i=0}^{n} \mu^{i} x_{i}(t), \quad t \in \bar{\mathcal{T}}.
\end{equation}

\subsection{Proof of Theorem \ref{MainThm}}

To prove Theorem \ref{MainThm} and estimate its accuracy \eqref{NorderEstim1}–\eqref{NorderEstim3}, we use the asymptotic method of inequalities \cite{b7}. First, we recall the definition of upper and lower solutions and their role in the construction of solution  \eqref{equat1} \cite{b7,b8,NEFEDOV201390}.

\begin{definition} \label{Lemma1}
 The functions  $   \beta (x,t,\mu) $ and  $ \alpha (x,t, \mu)$ are called upper and lower solutions of problem \eqref{equat1} if they are continuous, twice continuously differentiable in $x$, continuously differentiable in $t$, and for a sufficiently small $\mu$, satisfy the following conditions:

  \begin{itemize}[leftmargin=1cm]
 \item[\textbf{(C1):}] $\alpha(x,t,\mu)\leq \beta(x,t,\mu)$ for $(x,t)\in \bar{\Omega}\times \bar{\mathcal{T}}.$
  \item[\textbf{(C2):}]  $  \displaystyle L[ \alpha]:=\mu\frac{\partial^{2}\alpha}{dx^{2}}-\frac{\partial\alpha}{\partial t}+k\alpha \frac{\partial\alpha}{\partial x}-f(x)\geq 0, \quad (x,t)\in \bar{\Omega}\times \bar{\mathcal{T}};$
   \item[\qquad \ \ ] $ \displaystyle L[\beta]:=\mu\frac{\partial^{2}\beta}{dx^{2}}-\frac{\partial\beta}{\partial t}+k \beta \frac{\partial\beta}{\partial x}-f(x )\leq 0, \quad (x,t)\in \bar{\Omega}\times \bar{\mathcal{T}}.$
 \item[\textbf{(C3):}] $\alpha (0,t,\mu)\leq u^{l}\leq \beta(0,t,\mu), \quad \alpha(1,t,\mu)\leq u^{r}\leq \beta(1,t,\mu)$.
\end{itemize}
\end{definition}

\begin{lemma} \label{Lemma2}
(\cite{b8}) Let there be an upper $   \beta (x,t,\mu) $ and a lower $ \alpha (x,t, \mu)$ solution to problem \eqref{equat1} satisfying conditions (C1)–(C3) in Definition \ref{Lemma1}. Then, under Assumptions 1–4, there exists a solution $u(x,t,\mu)$ to problem \eqref{equat1} that satisfies the inequalities
$$
\alpha(x,t,\mu)\leq u(x,t,\mu)\leq \beta(x,t,\mu),\ (x,t)\in \bar{\Omega}\times \bar{\mathcal{T}}.
$$
Moreover, the functions $\beta(x, t, \mu)$ and $\alpha (x, t, \mu)$ satisfy the following estimates:
\begin{align}
& \beta(x,t,\mu)-\alpha(x,t,\mu) =\mathcal{O}(\mu^{n}), \label{BetaMinusAlpha} \\
& u(x,t,\mu)=\alpha(x,t,\mu)+\mathcal{O}(\mu^{n})=U_{n-1}(x,t,\mu)+\mathcal{O}(\mu^{n}). \label{UMinusAlpha}
\end{align}
 \end{lemma}

\begin{lemma} \label{Lemma3}
(\cite{b7,NEFEDOV201390}) Lemma \ref{Lemma2} also remains valid in the case in which the functions $\alpha(x,t,\mu)$ and $\beta(x,t,\mu)$ are continuous and their derivatives with respect to $x$ have discontinuities  from the class $C^2$ on some curve $x_{t.p}$, and the limit values of the derivatives on the curve $x_{t.p}$ satisfy the following conditions:
\begin{itemize}[leftmargin=1cm]
\item[\textbf{(C4):}] $\displaystyle \frac{\partial\alpha^{l}}{\partial x} \vert_{x=\underline{x}(t,\mu)}-\frac{\partial\alpha^{r}}{\partial x} \vert_{x=\underline{x}(t,\mu)}\leq 0, \quad
 \frac{\partial\beta^{l}}{\partial x} \vert_{x=\overline{x}(t,\mu)}-\frac{\partial\beta^{r}}{\partial x} \vert_{x=\overline{x}(t,\mu)}\geq 0.$\\
\end{itemize}
\end{lemma}

The proofs of Lemmas \ref{Lemma2}–\ref{Lemma3} can be found in \cite{b7,b8}. Thus, to prove Theorem \ref{MainThm}, it is necessary to construct the lower and upper solutions  $ \alpha (x,t, \mu)$  and $   \beta (x,t,\mu) $. Under conditions (C1)–(C4) for $ \alpha (x,t, \mu)$ and $   \beta (x,t,\mu) $, estimates \eqref{NorderEstim1}, \eqref{NorderEstim2} will follow directly from Lemma \ref{Lemma2}. Estimate \eqref{NorderEstim3} can be obtained by solving equation \eqref{proofzn2.20} for $z_{n}(x, t, \mu):= u(x, t, \mu)-U_{n}(x, t, \mu)$ using Green's function.

Now, we begin the proof of Theorem \ref{MainThm}.

\begin{proof}
Following the idea in \cite{b6},
we construct the upper and lower solutions $\alpha^{l}$, $\alpha^{r}$, $\beta^{l}$, $\beta^{r}$ and curves $\overline{x}$,  $\underline{x}$ as a modification of asymptotic representation \eqref{asymptoticnorder}.

We introduce a positive function $ \rho (t) $, which will be defined later in \eqref{phoequat}, and use the notations  $ \underline{\rho} (t)=\rho (t) $ and $ \overline{\rho} (t)=-\rho (t) $ to enable us to define the curves $\overline{x}(t)$ and $\underline{x}(t)$ in the form
\begin{align}\label{curveX}
\displaystyle \overline{x}(t,\mu )=\sum_{i=0}^{n+1} \mu^i x_{i}(t)+\mu^{n+1} \overline{\rho} (t),  \quad
\displaystyle \underline{x}(t,\mu )= \sum_{i=0}^{n+1} \mu^i x_{i}(t)+\mu^{n+1}\underline{\rho} (t).
\end{align}

Then,
\begin{equation} \label{UpperLowerV}
\displaystyle \overline{v}(t)=\frac{d\overline{x}}{dt},  \quad  \underline{v}(t)=\frac{d\underline{x}}{dt}.
\end{equation}

We introduce the stretched variables
\begin{equation} \label{UpperLowerXi}
\overline{\xi}=\frac{x- \overline{x}(t,\mu)}{\mu}, \quad \underline{\xi}=\frac{x- \underline{x}(t,\mu)}{\mu}.
\end{equation}

The upper and lower solutions of problem \eqref{equat1} will be constructed separately in the domains $ \bar{K}^{l}, \bar{K}^{r}$ and $\bar{M}^{l},\bar{M}^{r}$, in which the curves $\overline{x} (t)$ and $\underline{x} (t)$ divide the domain $ \bar{\Omega}\times \bar{\mathcal{T}} $:
\begin{align}
\label{beta}
\beta(x, t, \mu)= \begin{cases}
\beta^{l}(x,t,\mu ), \ (x,t)\in \bar{K}^{l}:=\{(x,t): x\in [0, \overline{x}(t,\mu)],  t \in \bar{\mathcal{T}} \},  \\
\beta^{r}(x,t,\mu ), \ (x,t)\in \bar{K}^{r}:=\{(x,t): x\in [\overline{x}(t,\mu), 1],  t \in \bar{\mathcal{T}} \}, 
\end{cases}
\end{align}
\begin{align}
\label{alpha}
\alpha(x, t, \mu)= \begin{cases}
\alpha^{l}(x,t,\mu ), \quad (x,t)\in \bar{M}^{l}:=\{(x,t): x\in [0, \underline{x}(t,\mu)],  t \in \bar{\mathcal{T}} \}, \\
\alpha^{r}(x,t,\mu ), \quad (x,t)\in \bar{M}^{r}:=\{(x,t): x\in [\underline{x}(t,\mu), 1] ,  t \in \bar{\mathcal{T}} \}.
\end{cases}
\end{align}

We will match the functions $\beta^{l}(x, t, \mu), \beta^{r}(x, t, \mu)$ and $\alpha^{l}(x, t, \mu), \alpha^{r}(x, t, \mu)$ on the curves $\overline{x}$ and $\underline{x}$, respectively, so that $\beta(x, t, \mu)$ and $\alpha(x, t, \mu)$ are continuous on these curves and the following equations hold:
\begin{align}
\begin{split} \label{sewingeq}
\displaystyle \beta^{l}(\overline{x}(t,\mu), t, \mu)=\beta^{r}(\overline{x}(t,\mu), t, \mu)=\frac{\varphi^{l}(\overline{x}(t,\mu))+\varphi^{r}(\overline{x}(t,\mu))}{2}, \\
\displaystyle \alpha^{l}(\underline{x}(t,\mu), t, \mu)=\alpha^{r}(\underline{x}(t,\mu), t, \mu)=\frac{\varphi^{l}(\underline{x}(t,\mu))+\varphi^{r}(\underline{x}(t,\mu))}{2}.
\end{split}
\end{align}

Note that we do not match the derivatives of the upper and lower solutions on the curves $\overline{x}(t)$ and $\underline{x}(t)$, and so the derivatives $\partial \beta / \partial x$  and $\partial \alpha / \partial x $ have discontinuity points, and therefore we need condition (C4) to hold.

We construct the functions $\beta^{l,r}$ and $\alpha^{l,r}$ in the following forms:
\begin{align}
\begin{split} \label{beta2}
\beta^{l,r}= U_{n+1}^{l,r} \vert_{\overline{\xi},\overline{x}}+\mu^{n+1} \left(\epsilon^{l,r}(x)+q_{0}^{l,r}(\overline{\xi}, t)+\mu q_{1}^{l,r}(\overline{\xi},t) \right), \\
\alpha^{l,r}=U_{n+1}^{l,r} \vert_{\underline{\xi},\underline{x}}-\mu^{n+1} \left( \epsilon^{l,r}(x)+q_{0}^{l,r}(\underline{\xi},t)+\mu q_{1}^{l,r}(\underline{\xi}, t)\right) ,
\end{split}
\end{align}
where the functions $\epsilon^{l,r}(x)$ should be designed in such a way that the condition (C2) is satisfied for $\beta^{l,r}$ and $\alpha^{l,r}$ in \eqref{beta2}. The functions $q_{0}^{l,r}(\overline{\xi},t) $ eliminate residuals of  order $\mu^n$ arising in $L[\beta]$ and $L[\alpha]$ and residuals of  order $\mu^{n+1}$ under the condition of continuous matching of the upper solution \eqref{sewingeq}, which arise as a result of modifying the regular part by adding  $\epsilon^{l,r}(x)$. The functions $q_{1}^{l,r}(\overline{\xi}, t) $ eliminate residuals of  order $\mu^{n+1}$ arising in $L[\beta]$  as we add  $\epsilon^{l,r}(x)$ and $q_{0}^{l,r} (\overline{\xi}, t) $.

Now, we define the functions $\epsilon^{l,r}(x)$ from the following equations:
\begin{align} \label{epsiloneq}
\begin{split}
k \frac{d \epsilon^{l,r} (x)}{d x} \varphi^{l,r} (x) +\epsilon^{l,r} (x) k \frac{d \varphi^{l,r} (x)}{d x}= -R  ,\\
\epsilon^{l}(0)=R^{l}, \quad \epsilon^{r} (1)= R^{r},
\end{split}
\end{align}
where $R, R^{l}, R^{r}$ are some positive values, which will be determined later. The functions $\epsilon^{l,r} (x)$ can be determined explicitly:
\begin{align} \label{epsiloneqexplicit}
\begin{split}
&\epsilon^{l} (x)= \frac{1}{\varphi^{l} (x)} (R^l \varphi^{l}(0)-\frac{Rx}{k}),\\
&\epsilon^{r} (x)= \frac{1}{ \varphi^{r} (x)} ( R^r \varphi^{r}(1)+ \frac{R}{k}(1-x)).
\end{split}
\end{align}

Since $\varphi^{l} (x) <0 $ and $\varphi^{r} (x) >0 $, $\epsilon^{l,r} (x) > 0 $ for $x \in \bar{\Omega}$.

 We define the functions $q_{0}^{l,r}(\overline{\xi},t) $ as solutions of the equations
\begin{align} \label{q0}
\begin{split}
\frac{\partial^2 q_{0}^{l,r}}{\partial \xi^2}+ \left( v_0(t) + k(\varphi^{l,r} (x_0(t)) + Q_{0}^{l,r}(\overline{\xi}, t)) \right)\frac{\partial q_{0}^{l,r}}{\partial \xi}+k q_{0}^{l,r} \Upsilon^{l,r}(\overline{\xi}, t) = H_{q0}^{l,r}(\overline{\xi}, t), \\
\end{split}
\end{align}
where
\begin{equation} \label{fq0}
H_{q0}^{l,r}(\overline{\xi}, t)=   \Upsilon^{l,r}(\overline{\xi}, t) \left( k \left( -\overline{\rho}(t) \frac{d \varphi^{l,r}}{d x}(x_0 (t)) - \epsilon^{l,r}(x_0 (t)) \right)-  \frac{d \overline{\rho}(t)}{dt} \right) .
\end{equation}

The boundary conditions for $q_{0}^{l,r}(\overline{\xi},t)$ follow from the conditions of continuous matching of the upper solution \eqref{sewingeq}, with the following conditions in  $\overline{\xi} = 0$ for functions $Q_{i}^{l,r}(\overline{\xi}, t) $:
\begin{align}  \label{q0conditions}
 q_{0}^{l,r}(0,t)=-\epsilon^{l,r}(x_0 (t))- \displaystyle \overline{\rho} (t) \frac{d \varphi^{l,r}}{d x} (x_0 (t))\equiv p_{2}^{l,r} (t), \quad q_{0}^{l}(-\infty,t) = 0, \quad q_{0}^{r}(+\infty,t) = 0 .
\end{align}

We can write the functions $q_{0}^{l}(\overline{\xi}, t)$ in this explicit form:
\begin{equation} \label{q0explisit}
q_{0}^{l,r} (\overline{\xi}, t)=z^{l,r}(\overline{\xi},t) \left( p_{2}^{l,r} (t)   -\int_{0}^{\overline{\xi}} \frac{1}{z^{l,r}(s,t)} \int_{s}^{\mp \infty} H_{q0}^{l,r} (\eta,t) d\eta ds \right).
\end{equation}

 We define the functions $q_{1}^{l,r}(\overline{\xi}, t) $ from the following equations:
\begin{align} \label{q1}
\begin{split}
& \frac { \partial^{2} q_{1}^{l,r}}{ \partial {\overline{\xi}}^{2}} + \left( k (Q_{0}^{l,r} ( \overline{\xi},t ) +\varphi^{l,r}( x_{0} ( t ) )) +v_{0} ( t )  \right) \frac {\partial q_{1}^{l,r} }{\partial \overline{\xi}} +kq_{1}^{l,r}   \Upsilon^{l,r}(\overline{\xi}, t)  = \\ & \qquad
k \left( -\overline{\rho} ( t ) \frac{\partial Q_{1}^{l,r} }{\partial   \overline{\xi}} ( \overline{\xi},t ) + ( -\overline{\xi} -x_{1} ( t )  ) \frac{\partial  q_{0}^{l,r}}{\partial   \overline{\xi}}  ( \overline{\xi}, t ) -q_{0}^{l,r}( \overline{\xi}, t )  \right) \frac{d \varphi^{l,r}  }{d  x }   ( x_{0} ( t ) ) \\ 
& \qquad
-k\overline{\rho}( t )  \left(  ( \overline{\xi}+x_{1} ( t )  )  \Upsilon^{l,r}(\overline{\xi}, t)  +Q_{0}^{l,r} ( \overline{\xi},t )  \right)  \frac { d^{2}  \varphi^{l,r}}{ d {x}^{2}}( x_{0} ( t )  )  \\
 & \qquad
-k \left( \left( \overline{\xi}+x_{1} ( t )  \right) \frac{d \epsilon^{l,r} }{d   x}  ( x_{0} ( t ) )  + \overline{\rho} ( t ) \frac{d u_{1}}{d   x}  ( x_{0} ( t ) ) \right)  \Upsilon^{l,r}(\overline{\xi}, t) \\ 
& \qquad
- \left( k\epsilon^{l,r} ( x_{0} ( t )  ) +kq_{0}^{l,r}( \overline{\xi} )+ \frac{d \overline{\rho} ( t )}{d t} \right) \frac{\partial Q_{1}^{l,r} }{\partial  \overline{\xi} }  ( \overline{\xi},t )+\frac{\partial q_0}{\partial t} ( \overline{\xi},t ) \\ 
& \qquad
 -k\frac{d \epsilon^{l,r} }{d  x } ( x_{0} ( t )  ) Q_{0}^{l,r} ( \overline{\xi},t ) -\frac{\partial  q_{0}^{l,r}}{\partial  \overline{\xi} }  ( \overline{\xi}, t )  \left( k(u_{1} ( x_{0} ( t ) ,t )+Q_{1}^{l,r} ( \overline{\xi},t )) +v_{1}  ( t )  \right),
\end{split}
\end{align}
with the boundary conditions
\begin{align*}
q_{1}^{l,r}(0,t)=0, \ q_{1}^{l,r}(\overline{\xi},t) \rightarrow 0 \ \text{for} \  \overline{\xi} \rightarrow \mp \infty.
\end{align*}

Replacing $\overline{\rho}$ with $\underline{\rho}$ and $\overline{\xi}$ with $\underline{\xi}$ in \eqref{q0}–\eqref{q1}, we define the functions $q_{0}^{l,r}(\underline{\xi}, t)$ and $q_{1}^{l,r}(\underline{\xi}, t)$ that appear in the functions $\alpha^{l,r}$.

The functions $q_{0}^{l,r}$ and $q_{1}^{l,r}$ satisfy exponential estimates of type \eqref{equat22} and \eqref{equat23}.

Now, we need to show that the functions $\beta(x,t,\mu)$ and $\alpha(x,t,\mu)$ are upper and lower solutions to problem \eqref{equat1}. To do this, we check all conditions (C1)–(C4).

Condition (C1) is checked in the same way as in \cite{b7}. Using equations \eqref{beta}, \eqref{alpha}, and \eqref{beta2}, it is possible to verify that $\beta-\alpha>0$ for each of the regions: $[0, \overline{x}(t,\mu)], [\overline{x}(t,\mu), \underline{x}(t,\mu)], [\underline{x}(t,\mu),1 ]  $. 

The method of constructing the upper and lower solutions implies the following inequalities:
\begin{equation*}
L[\beta]=-\mu^{n+1} R + \mathcal{O}(\mu^{n+2})<0, \quad L[\alpha]=\mu^{n+1} R + \mathcal{O}(\mu^{n+2})>0,
\end{equation*}
where $R$ is a constant from \eqref{epsiloneq}. This verifies condition (C2).

Condition (C3) is satisfied for sufficiently large values $R^l$ and $R^r$ in the boundary conditions of equation \eqref{epsiloneq}.

We now check condition (C4) for the upper solutions $\beta^{l,r}$. Because of the matching conditions \eqref{matching}, \eqref{matchingfirstord} (and up to order $n+1$), the coefficients for $\mu^{i}$ ($i = 1,\cdots,n$) are equal to zero, and the coefficient at $\mu^{n + 1}$ includes only the terms resulting from the modification of the asymptotics:
\begin{align}
\begin{split}
&\mu \left( \frac{\partial\beta^{l}}{\partial x}-\frac{\partial\beta^{r}}{\partial x} \right)\Big \vert_{x=\overline{x}(t)} = \mu^{n+1} \left( \frac{\partial {q_0}^l}{\partial \xi}(0,t) - \frac{\partial q_{0}^r}{\partial \xi}(0,t) \right) +\mathcal{O}(\mu^{n+2})\\
\end{split}.
\end{align}

Using the explicit solution for $  {q_0}^{l,r}(0,t)$ \eqref{q0explisit}, we find
\begin{align}
\begin{split}
& \frac{\partial {q_0}^l}{\partial \xi}(0,t) - \frac{\partial q_{0}^r}{\partial \xi}(0,t)  =(v_0(t)+k \varphi^l(x_0(t))) \left( \epsilon^{l}(x_0 (t))-\displaystyle \rho (t) \frac{d \varphi^{l}}{d x} (x_0 (t))\right) - \frac{d \rho (t)}{dt} \varphi^{l} (x_0 (t))  \\ &\qquad\qquad
- (v_0(t)+k \varphi^r(x_0(t))) \left( \epsilon^{r}(x_0 (t))- \displaystyle \rho (t) \frac{d \varphi^{r}}{d x} (x_0 (t))\right) + \frac{d \rho (t)}{dt} \varphi^{r} (x_0 (t)) \\ &\qquad
=\frac{ d \rho(t)}{d t}   (\varphi^{r}(x_0 (t)) - \varphi^{l}(x_0 (t)))-\frac{k}{2} \left( \varphi^{r}(x_0 (t))- \varphi^{l}(x_0 (t))\right) \left( \epsilon^l +\epsilon^r \right) \\ &\qquad\qquad
+ \rho(t) \left( \frac{k}{2}\left(\varphi^{r}(x_0(t))-\varphi^{l}(x_0(t))  \right) \left( \frac{d \varphi^{r} }{d x}(x_0(t))+ \frac{d \varphi^{l} }{d x}(x_0(t)) \right) \right)  .
\end{split}
\end{align}

We choose the function $\rho (t) $ as a solution to the problem
\begin{align} \label{phoequat}
\frac{d \rho (t)}{dt} (\varphi^{r}(x_0 (t))- \varphi^{l}(x_0 (t)) ) = -\Phi_1 (t) \rho(t)+F(t)+\sigma, \quad \rho (0)=\rho^0, \quad t \in \bar{\mathcal{T}},
\end{align}
where $F(t)=  \frac{k}{2} \left( \varphi^{r}(x_0 (t))- \varphi^{l}(x_0 (t))\right) \left( \epsilon^l +\epsilon^r \right)$. Since the function $F(t)$ and the constants $\sigma$ and $\rho^0$ are positive, the solution  $\rho (t) $ to equation \eqref{phoequat} is also positive.

For such $\rho (t)$, we obtain: 
\begin{equation}
\mu \left( \frac{\partial\beta^{l}}{\partial x}-\frac{\partial\beta^{r}}{\partial x} \right)\Big \vert_{x=\overline{x}(t)}=\mu^{n+1} \sigma +\mathcal{O}(\mu^{n+2}) >0.
\end{equation}

Similarly, condition (C4) is satisfied for the functions $\alpha^{l,r}$, and the constructed upper and lower solutions guarantee the existence of a solution $u (x, t, \mu)$ to problem \eqref{equat1}, satisfying the inequalities
\begin{equation}
 \alpha(x,t,\mu)\leq u(x,t, \mu)\leq \beta(x,t,\mu).
\end{equation}
In addition, estimates \eqref{NorderEstim1},\eqref{NorderEstim2} are valid.

We now show that estimate \eqref{NorderEstim3}  also holds. To do this, we estimate the difference  $z_{n}(x, t, \mu)\equiv u(x, t, \mu)-U_{n}(x, t, \mu) $; the function $z_{n}(x, t, \mu)$ satisfies the equation
\begin{align} \label{proofzn2.20}
\begin{split}
\mu \frac{\partial^{2}z_{n}}{\partial x^{2}}-\frac{\partial z_{n}}{\partial t}- \left( k U_{n} \frac{\partial U_{n}}{\partial x} - k u \frac{\partial u}{\partial x}  \right)  =\mu^{n+1}\psi(x, t, \mu)
\end{split}
\end{align}
for $(x, t)\in \bar{\Omega}\times \bar{\mathcal{T}}$, with zero boundary conditions, where $\lvert \psi(x, t, \mu) \rvert \leq c_{1}$. Using the estimates from Lemma \ref{Lemma2}, we obtain
\begin{equation} \label{proof2.21}
z_{n}(x, t, \mu)= u(x, t, \mu)-U_{n}(x, t, \mu)\leq O (\mu^{n+1}).
\end{equation}

The second term of equation \eqref{proofzn2.20} can be represented in the form
\begin{equation}
 k U_{n} \frac{\partial U_{n}}{\partial x}- k u \frac{\partial u}{\partial x}=\frac{\partial}{\partial x}\int_{u}^{U_{n}} (k s) ds.
\end{equation}

We rewrite \eqref{proofzn2.20} in the following form:
\begin{align} \label{proof2.22}
\begin{split}
\frac{\partial^{2}z_{n}}{\partial x^{2}}- \frac{1}{\mu} \frac{\partial z_{n}}{\partial t}-K z_{n}=-K z_{n}+\frac{1}{\mu}\frac{\partial}{\partial x}\int_{u}^{U_{n}} (k s) ds + \mu^{n}\psi(x, t, \mu).
\end{split}
\end{align}

We define
$$
r(x, t, \mu):=\mu^{n}\psi(x, t, \mu),
$$
and, changing the variable to $\tilde{t}=\mu t$,  we can rewrite \eqref{proof2.22} in the following form:
\begin{align} \label{proof2.22v2}
\begin{split}
\frac{\partial^{2}z_{n}}{\partial x^{2}}- \frac{\partial z_{n}}{\partial \tilde{t}}-K z_{n}=-K z_{n}+\frac{1}{\mu}\frac{\partial}{\partial x}\int_{u}^{U_{n}} (k s) ds
+ r(x, \frac{\tilde{t}}{\mu}, \mu).
\end{split}
\end{align}

Using a Green's function for the parabolic operator on the left-hand side of \eqref{proof2.22v2}, for any $ (x,t)\in \bar{\Omega}\times \bar{\mathcal{T}}$, $t_{0} \in [0,t)$, and $ (\zeta,\tau) \in \bar{\Omega} \times [0,\mu t)$ we obtain the representation for $z_{n}$ \cite{pao1992}:
\begin{align} \label{proof2.23}
\begin{split}
z_{n}=\int_{0}^{1}G(x, \mu t, \zeta, \mu t_{0})z_{n}(\zeta, \mu t_{0})d\zeta-\int_{\mu t_{0}}^{\mu t}d\tau\int_{0}^{1}G(x, \mu t, \zeta, \frac{\tau}{\mu} ) 
\displaystyle  \left(-K z_{n}(\zeta, \frac{\tau}{\mu})+r(\zeta, \frac{\tau}{\mu}, \mu)+\frac{1}{\mu}\frac{\partial}{\partial\zeta}\int_{u(\zeta,\frac{\tau}{\mu},\mu)}^{U_{n}(\zeta,\frac{\tau}{\mu},\mu)}(k s)ds \right)d\zeta.
\end{split}
\end{align}

Using integration by parts and the boundary conditions for $G$, we can transform the last term in \eqref{proof2.23} as follows:
\begin{align} \label{proof2.24}
\begin{split}
& \int_{\mu t_{0}}^{\mu t}d\tau\int_{0}^{1}G(x,\mu t, \zeta, \frac{\tau}{\mu})\frac{1}{\mu}\frac{\partial}{\partial\zeta}\int_{u(\zeta,\frac{\tau}{\mu},\mu)}^{U_{n}(\zeta,\frac{\tau}{\mu},\mu)}(k s)dsd\zeta  =-\int_{\mu t_{0}}^{\mu t}d\tau\int_{0}^{1}G_{\zeta}(x, \mu t, \zeta, \frac{\tau}{\mu})\frac{1}{\mu}\int_{u(\zeta,\frac{\tau}{\mu},\mu)}^{U_{n}(\zeta,\frac{\tau}{\mu},\mu)}(k s)dsd\zeta \\ & \qquad
=-\int_{\mu t_{0}}^{\mu t}d\tau\int_{0}^{1}G_{x}(x, \mu t, \zeta, \frac{\tau}{\mu})\frac{1}{\mu}\int_{u(\zeta,\frac{\tau}{\mu},\mu)}^{U_{n}(\zeta,\frac{\tau}{\mu},\mu)}(k s)dsd\zeta 
=-\displaystyle \frac{\partial}{\partial x} \left( \int_{\mu t_{0}}^{\mu t}d\tau\int_{0}^{1}G_{x}(x, \mu t, \zeta, \frac{\tau}{\mu})\frac{1}{\mu}\int_{u(\zeta,\frac{\tau}{\mu},\mu)}^{U_{n}(\zeta,\frac{\tau}{\mu},\mu)}(k s)dsd\zeta \right) .
\end{split}
\end{align}

Using \eqref{proof2.24}, we obtain from \eqref{proof2.23} the following representation for the derivative $\displaystyle \frac{\partial z_{n}}{\partial x}$:
\begin{multline} \label{proof2.25}
 \frac{\partial z_{n}}{\partial x}=\int_{0}^{1}G_{x}(x, \mu t, \zeta, \mu t_{0})z_{n}(\zeta, \mu t_{0})d\zeta 
-\int_{\mu t_{0}}^{\mu t}d\tau\int_{0}^{1}G_{x}(x, \mu t, \zeta, \frac{\tau}{\mu})\left(-K z_{n}(\zeta, \frac{\tau}{\mu})+r(\zeta, \frac{\tau}{\mu}, \mu)\right)d\zeta \\ 
+\displaystyle \frac{\partial^{2}}{\partial x^{2}} \left(\int_{\mu t_{0}}^{\mu t}d\tau\int_{0}^{1}G(x, \mu t, \zeta, \frac{\tau}{\mu})\frac{1}{\mu}\int_{u(\zeta,\frac{\tau}{\mu},\mu)}^{U_{n}(\zeta,\frac{\tau}{\mu},\mu)}(k s)dsd\zeta \right) .
\end{multline}

The validity of representation \eqref{proof2.25} follows from the estimates
\begin{equation*}\label{estimat1}
\left \lvert \int_{0}^{1}G_{x}(x, \mu t, \zeta, \mu t_{0})d\zeta \right \rvert \leq C, \quad \left \lvert \int_{\mu t_{0}}^{\mu t}d\tau\int_{0}^{1}G_{x}(x, \mu t, \zeta, \frac{\tau}{\mu} )d\zeta \right \rvert \leq C
\end{equation*}
and
\begin{equation*}\label{estimat2}
\left \lvert \frac{\partial^{2}}{\partial x^{2}}\int_{\mu t_{0}}^{\mu t}d\tau\int_{0}^{1}G(x, \mu t, \zeta, \frac{\tau}{\mu})\frac{1}{\mu}d\zeta \right \rvert \leq C,
\end{equation*}
which can be found, for example, in \cite[Page 49]{pao1992}. We find that the first and second terms of representation \eqref{proof2.25} have estimates $\mathcal{O}(\mu^{n+1})$ and $\mathcal{O}(\mu^{n})$, respectively. We also find that the last term in representation \eqref{proof2.25} can be estimated by
$$
\frac{1}{\mu} \left \lvert \int_{u(x,\mu t,\mu)}^{U_{n}(x,\mu t,\mu)}(k s)ds\right \rvert \leq \mathcal{O}(\mu^{n}) .
$$
Using these estimates, from \eqref{proof2.25} we obtain $\displaystyle \frac{\partial z_{n}}{\partial x}(x, \mu t, \mu)=\mathcal{O}(\mu^{n})$ for $(x, t)\in \bar{\Omega}\times \bar{\mathcal{T}}$. This completes the proof of Theorem \ref{MainThm}.  
\end{proof}

\subsection{Proof of Lemma \ref{LemmaC}}
By the assumptions of the lemma, we deduce that:
\begin{align}
\label{estimatel}
\left\| \varphi^{l}(x) - u (x,t) \right\|_{L^p(\bar{\Omega}^{l})}
& =  \left\| Q_0^{l}+\sum_{i=1}^{\infty} \mu^i \left( \bar{u}_i^{l}+Q_i^{l} \right) \right\|_{L^p(\bar{\Omega}^{l})}
\leq \left\| Q_0^{l} \right\|_{L^p(\bar{\Omega}^{l})} + \left\|\sum_{i=1}^{\infty} \mu^i \left( \bar{u}_i^{l}+Q_i^{l} \right) \right\|_{L^p(\bar{\Omega}^{l})} \nonumber \\ 
& \leq \left\| Q_0^{l} \right\|_{L^p(\bar{\Omega}^{l})} + 2 \mu \left\|\bar{u}_1^{l}+Q_1^{l} \right\|_{L^p(\bar{\Omega}^{l})}
< \mu ( 1 + 2\left\|\bar{u}_1^{l}+Q_1^{l} \right\|_{L^p(\bar{\Omega}^{l})}),
\end{align}
taking into account the bounds for considered small $\mu$ and the inequality $\lvert Q_{0}^{l} \rvert \leq \mu^2$ in the region $(0,x_0 - \Delta x/2)$. In the same way, we obtain:
\begin{align}
\label{estimatelD}
& \left\|  \frac{d \varphi^{l}(x)}{dx} - \frac{\partial u(x,t)}{\partial x} \right\|_{L^p(0,x_0 - \Delta x/2)}
=  \left\| \frac{\partial Q_0^{l}}{\partial x} +\sum_{i=1}^{\infty} \mu^i \frac{\partial }{\partial x} (\bar{u}_i^{l}+Q_i^{l}) \right\|_{L^p(0,x_0 - \Delta x/2)}
\leq \left\| \frac{\partial Q_0^{l}}{\partial x} \right\|_{L^p(0,x_0 - \Delta x/2)} + \left\|\sum_{i=1}^{\infty} \mu^i \frac{\partial }{\partial x} (\bar{u}_i^{l}+Q_i^{l}) \right\|_{L^p(0,x_0 - \Delta x/2)} \nonumber \\ 
& \quad  \leq  \frac{1}{\mu} \left\| \frac{\partial Q_0^{l}}{\partial \xi} \right\|_{L^p(0,x_0 - \Delta x/2)} + 2 \mu \left\| \frac{\partial }{\partial x} (\bar{u}_1^{l}+Q_1^{l}) \right\|_{L^p(0,x_0 - \Delta x/2)}
< \mu ( k \left\|P^{l}\right\|_{L^p(0,x_0 - \Delta x/2)} + 2\left\| \frac{\partial \bar{u}_{1}^{l}}{\partial x}\right\|_{L^p(0,x_0 - \Delta x/2)} +\frac{2}{\mu}\left\|\frac{\partial Q_1^{l}}{\partial \xi}  \right\|_{L^p(0,x_0 - \Delta x/2)}) \nonumber\\
& \quad \leq \mu ( 3k \left\|P^{l}\right\|_{L^p(0,x_0 - \Delta x/2)} + 2\left\| \frac{\partial \bar{u}_{1}^{l}}{\partial x}\right\|_{L^p(0,x_0 - \Delta x/2)} ),
\end{align}
where we used the estimates $\mu \left\| \frac{\partial Q_{1}^{l}}{\partial \xi} \right\| \leq \left\| \frac{\partial Q_{0}^{l}}{\partial \xi} \right\| $,  $ \lvert \frac{\partial Q_{0}^{l}}{\partial \xi} \rvert \leq -k P^{l}\mu^2- k\mu^4/2 < -k P^{l}\mu^2$  and $P^{l}(x_0(t))<0$ is defined in the equation \eqref{zeroordertransitionfunc}. By combining \eqref{estimatel} and \eqref{estimatelD}, we conclude that 
\begin{align*}
\left\| \varphi^{l}(x) - u (x,t) \right\|_{W^{1,p}(0,x_0 - \Delta x/2)} \leq C' \mu 
\end{align*}
with $C' :=  1 + 2\left\|\bar{u}_1^{l}+Q_1^{l} \right\|_{L^p(\bar{\Omega}^{l})} + 3k \left\|P^{l}\right\|_{L^p(0,x_0 - \Delta x/2)} + 2\left\| \frac{\partial \bar{u}_{1}^{l}}{\partial x}  \right\|_{L^p(0,x_0 - \Delta x/2)}$.

Similarly, we can derive the estimates in the right region $\bar{\Omega}^{r}$:

\begin{align}
\label{estimateR}
\left\| \varphi^{r}(x) - u (x,t) \right\|_{L^p(\bar{\Omega}^{r})}
& =  \left\| Q_0^{r}+\sum_{i=1}^{\infty} \mu^i \left( \bar{u}_i^{r}+Q_i^{r} \right) \right\|_{L^p(\bar{\Omega}^{r})}
< \mu ( 1 + 2\left\|\bar{u}_1^{r}+Q_1^{r} \right\|_{L^p(\bar{\Omega}^{r})}),
\end{align}

\begin{align}
\label{estimatelDR}
& \left\|  \frac{d \varphi^{r}(x)}{dx} - \frac{\partial u(x,t)}{\partial x} \right\|_{L^p(x_0 + \Delta x/2,1)}
=  \left\| \frac{\partial Q_0^{r}}{\partial x} +\sum_{i=1}^{\infty} \mu^i \frac{\partial }{\partial x} (\bar{u}_i^{r}+Q_i^{r}) \right\|_{L^p(x_0 + \Delta x/2,1)} \leq \mu ( 3k \left\|P^{r}\right\|_{L^p(x_0 + \Delta x/2,1)} + 2\left\| \frac{\partial \bar{u}_{1}^{r}}{\partial x}\right\|_{L^p(x_0 + \Delta x/2,1)} ).
\end{align}

By combining \eqref{estimateR} and \eqref{estimatelDR} we obtain:

\begin{align*}
\left\| \varphi^{r}(x) - u (x,t) \right\|_{W^{1,p}(x_0 + \Delta x/2,1)} \leq C'' \mu 
\end{align*}
with $C'' :=  1 + 2\left\|\bar{u}_1^{r}+Q_1^{r} \right\|_{L^p(\bar{\Omega}^{r})} + 3k \left\|P^{r}\right\|_{L^p(x_0 + \Delta x/2,1)} + 2\left\| \frac{\partial \bar{u}_{1}^{r}}{\partial x}  \right\|_{L^p(x_0 + \Delta x/2,1)}$, and the constant required for the lemma can be obtained by $C=C'+C''$.

\subsection{Proof of Proposition \ref{ProAsympErr} }

\begin{proof}
First, we note that the exact source function $f^*$ has the following representation according to equation \eqref{LeftRightSolution}:
\begin{align} \label{equatforf*}
f^*= \begin{cases}
\displaystyle k \varphi^{l} (x) \frac{d\varphi^{l} (x)}{dx}, \quad x \in (0,x_0(t) - \Delta x/2), \\ \\
\displaystyle k \varphi^{r} (x) \frac{d\varphi^{r} (x)}{dx}, \quad x \in (x_{0}(t)+\Delta x /2, 1).
  \end{cases}
\end{align}

Let $\Omega'=(0,x_0 - \Delta x/2)$. By using Lemma \ref{LemmaC}, we have
\begin{equation} \label{eq004IPPf}
\left\| \varphi^{l}(x) - u (x,t) \right\|_{L^p(\Omega')} \leq C \mu,
\end{equation}

\begin{equation}\label{0orderEstim1IPPf1}
\left\| \frac{d \varphi^{l}(x)}{dx} - \frac{\partial u(x,t)}{\partial x} \right\|_{L^p(\Omega')} \leq C \mu ,
\end{equation}

\begin{equation}\label{0orderEstim1IPPf2}
\left\| \frac{d \varphi^{l}(x)}{dx} \right\|_{L^p(\Omega')} \leq \left\| \frac{d \varphi^{l}(x)}{dx} - \frac{\partial u(x,t)}{\partial x} \right\|_{L^p(\Omega')} + \left\|\frac{\partial u(x,t)}{\partial x} \right\|_{L^p(\Omega')} \leq C \mu + \left\|\frac{\partial u(x,t)}{\partial x} \right\|_{L^p(\Omega')}.
\end{equation}

From estimates \eqref{eq004IPPf}–\eqref{0orderEstim1IPPf2}, we conclude that
\begin{align} \label{proofPro1}
\begin{split}
& \frac{1}{k} \left\|f^* - f_0 \right\|_{L^p(\Omega')} = \left\| \varphi^{l}(x) \frac{d \varphi^{l}(x)}{dx} - u (x,t) \frac{\partial u (x,t)}{\partial x} \right\|_{L^p(\Omega')} \\ & \quad \leq \left\| \varphi^{l}(x) \frac{d \varphi^{l}(x)}{dx}  - u (x,t) \frac{d \varphi^{l}(x)}{dx} \right\|_{L^p(\Omega')} +   \left\| u (x,t) \frac{d \varphi^{l}(x)}{dx}  - u (x,t) \frac{\partial u (x,t)}{\partial x} \right\|_{L^p(\Omega')}
 \\ & \quad \leq \left\| \varphi^{l}(x) - u (x,t) \right\|_{L^p(\Omega')} \left\| \frac{d \varphi^{l}(x)}{dx} \right\|_{L^p(\Omega')} +   \left\| u (x,t) \right\|_{L^p(\Omega')} \left\| \frac{d \varphi^{l}(x)}{dx}  -\frac{\partial u (x,t)}{\partial x} \right\|_{L^p(\Omega')}\\ & \quad
  \leq C\left( C+\left\|u(x,t) \right\|_{W^{1,p}(\Omega')} \right) \mu =: \frac{c_1}{k} \mu.
\end{split}
\end{align}

Following exactly the same lines, we also derive the inequality
\begin{equation} \label{proofPro2}
\left\|f^* - f_0 \right\|_{L^p(x_0 + \Delta x/2,1)}\leq c_2 \mu
\end{equation}
with a constant $c_2$. In addition, since $\Delta x\sim\mu \lvert \ln\mu \rvert$ (see also \eqref{DeltaX}), we have  
\begin{align} \label{proofPro3}
\displaystyle \left\|f^* - f_0 \right\|_{L^p(x_0 - \Delta x/2, x_0 + \Delta/2)}\leq
\left\|f^* \right\|_{L^p(x_0 - \Delta x/2, x_0 + \Delta x/2)} + \left\| f_0 \right\|_{L^p(x_0 - \Delta x/2, x_0 + \Delta x/2)} 
\leq c_3 \mu \lvert \ln \mu \rvert
\end{align}
with $c_3=\left\|f^* \right\|_{C(\Omega)} + \left\| f_0 \right\|_{C(\Omega)}$. By combining  \eqref{proofPro1}–\eqref{proofPro3}, we deduce that
\begin{align*}
\begin{split}
\left\|f^* - f_0 \right\|^p_{L^p(0,1)}  =  \left\|f^* - f_0 \right\|^p_{L^p(0,x_0 - \frac{\Delta x}{2})} + \left\|f^* - f_0 \right\|^p_{L^p(x_0 - \frac{\Delta x}{2},x_0 + \frac{\Delta x}{2})}  + \left\|f^* - f_0 \right\|^p_{L^p(x_0 + \frac{\Delta x}{2},1)}   \leq
c^p_1 \mu^p + c^p_2 \mu^p + c^p_3 \mu^p \lvert \ln \mu \rvert^p ,
\end{split}
\end{align*}
which yields required estimate \eqref{f0Ineq} with $C_1 = \left( c^p_1  + c^p_2 + c^p_3 \right)^{1/p}$.
\end{proof}

\subsection{Proof of Proposition \ref{NoisyErr}}

Without loss of generality, we assume that $u^\delta_0(t)\equiv u(x_0,t)$ and $u^\delta_n(t)\equiv u(x_n,t)$. Otherwise, we can consider the function
\begin{equation}
\label{uBar}
\bar{u}(x,t) = u(x,t) + u^\delta_0(t) - u(0,t) + b(t) x,
\end{equation}
where $b(t)= u^\delta_n(t) - u(1,t) + u(0,t) - u^\delta_0(t)$. It is clear that $\bar{u}(0,t)=u^\delta_0(t)$ and $\bar{u}(1,t)=u^\delta_n(t)$. All assertions below  hold according the triangle inequality
\begin{equation}
\label{triangleIneq}
\begin{array}{ll}
\|u^\varepsilon(\cdot,t)-u(\cdot,t)\|_{L^2(\Omega)} \leq  \|u^\varepsilon(\cdot,t)-\bar{u}(\cdot,t)\|_{L^2(\Omega)} + \|\bar{u}(\cdot,t)-u(\cdot,t)\|_{L^2(\Omega)} \\ \qquad = \|u^\varepsilon(\cdot,t)-\bar{u}(\cdot,t)\|_{L^2(\Omega)} + \lvert b(t)\rvert \leq \|u^\varepsilon(\cdot,t)-\bar{u}(\cdot,t)\|_{L^2(\Omega)} + 2\delta.
\end{array}
\end{equation}

\begin{proof}
Let $e(x,t):= u^\varepsilon(x,t) - u(x,t)$. From the definition of $u^\varepsilon(x,t)$ in \eqref{uAlpha}, we have $e(0,t)=e(1,t)=0$ for all $t\in \bar{\mathcal{T}}$. Consequently, from the Dirichlet–Poincare inequality, we obtain, for every $t\in \bar{\mathcal{T}}$,
\begin{equation}
\label{PoincareIneq}
\|e(x,t)\|^2_{L^2(\Omega)}  \leq 4 \left\|\frac{\partial e(x,t)}{\partial x }\right\|^2_{L^2(\Omega)}.
\end{equation}

For every $t\in \bar{\mathcal{T}}$, let $u^s(x,t)$ be the natural cubic spline over $\Theta$ that interpolates the exact data
$u(x,t)$ at the grid $\Theta$. Let $e_1(x,t)=u^\varepsilon(x,t) - u^s(x,t)$ and $e_2(x,t)=u^s(x,t)- u(x,t)$. It is clear that $e(x,t)=e_1(x,t)+e_2(x,t)$.

According to \cite[Lemmas 4.1, 4.2]{HankeScherzer2001}, for a fixed $t$, the following holds:
\begin{equation}
\label{Lemma2Ineq}
\left\|\frac{\partial e_2(x,t)}{\partial x }\right\|^2_{L^2(\Omega)} = \left\|\frac{\partial u^s(x,t)}{\partial x } - \frac{\partial u(x,t)}{\partial x } \right\|^2_{L^2(\Omega)}  \leq \frac{h}{\pi} \left\|\frac{\partial^2 u(x,t)}{\partial x^2 }\right\|^2_{L^2(\Omega)}.
\end{equation}
Moreover, for each $t$, $\frac{\partial^2 u^s(x,t)}{\partial x^2 }$ is the best approximation of $u(x,t)$ in $L^2(\Omega)$ from the space of linear splines over $\Theta$, i.e. the following identity holds:
\begin{equation}
\label{Lemma1Eq}
\left\|\frac{\partial^2 e_2(x,t)}{\partial x^2 } \right\|^2_{L^2(\Omega)} + \left\|\frac{\partial^2 u^s(x,t)}{\partial x^2 }\right\|^2_{L^2(\Omega)} = \left\|\frac{\partial^2 u(x,t)}{\partial x^2 }\right\|^2_{L^2(\Omega)}.
\end{equation}

On the other hand, for each $t\in \bar{\mathcal{T}}$, let $\chi(x,t)$ be the best approximating piecewise constant spline of $\frac{\partial e_1(x,t)}{\partial x }$ in $\mathcal{L}^2(\Omega)$, i.e.
\begin{equation} \label{chiequation}
\chi \vert_{\left(x_{i-1},x_{i}\right)} = \chi_{i} =
\frac{1}{h}\int_{x_{i-1}}^{x_{i}} \frac{\partial e_1}{\partial x } dx.
\end{equation}

Then, we obtain, together with $e_1(0,t)=e_1(1,t)=0$ for a.e. $t\in \bar{\mathcal{T}}$,
\begin{align} \label{normforderivativeofe}
\left\| \frac{\partial e_1}{\partial x } \right\|^2 &=\int_{0}^{1} \frac{\partial e_1}{\partial x }\left(\frac{\partial e_1}{\partial x }-\chi\right) dx + \int_{0}^{1} \frac{\partial e_1}{\partial x }
\chi dx \\
&= \int_{0}^{1} \frac{\partial e_1}{\partial x }\left(\frac{\partial e_1}{\partial x }- \chi\right) dx + \sum_{i=1}^{n} \chi_{i}
\int_{x_{i-1}}^{x_i} \frac{\partial e_1}{\partial x } dx  \nonumber
\\
&=\int_{0}^{1} \frac{\partial e_1}{\partial x }\left(\frac{\partial e_1}{\partial x }-\chi\right)dx +
\sum_{i=1}^{n}\chi_i\left(e_1\left(x_i,t\right)-e_1\left(x_{i-1},t\right)\right)\nonumber
\\
&=\int_{0}^{1} \frac{\partial e_1}{\partial x }\left(\frac{\partial e_1}{\partial x }-\chi\right)dx +
\sum_{i=1}^{n-1}e_1\left(x_i,t\right)\left(\chi_i-\chi_{i+1}\right)+e_1\left(1,t \right) \chi_n - e_1\left(0 ,t\right)\chi_1 \nonumber
\\
&= \int_{0}^{1} \frac{\partial e_1}{\partial x }\left(\frac{\partial e_1}{\partial x }-\chi\right)dx +
\sum_{i=1}^{n-1}e_1\left(x_i,t\right)\left(\chi_i-\chi_{i+1}\right)=:I_1+I_2.
\end{align}

From the approximation property of piecewise constant splines (cf. \cite[Theorem 6.1]{Schumaker1981}), we have
\begin{equation}
\left\Vert \frac{\partial e_1}{\partial x }-\chi \right\Vert_{L^2(\Omega)} \leq h \left\Vert \frac{\partial^2 e_1}{\partial x^2 } \right\Vert_{L^2(\Omega)} \nonumber,
\end{equation}
which implies, together with the Cauchy–Schwarz inequality, that
\begin{equation}
I_1 \leq \left\Vert \frac{\partial e_1}{\partial x } \right\Vert \left\Vert \frac{\partial e_1}{\partial x }-\chi \right\Vert \leq h \left\Vert \frac{\partial e_1}{\partial x } \right\Vert \,
\left\Vert  \frac{\partial^2 e_1}{\partial x^2 }  \right\Vert. \nonumber
\end{equation}

Since $u^\varepsilon(x,t)$ stands for a minimizer of \eqref{uAlpha}, we have \begin{align*}
\begin{array}{ll}
\displaystyle \delta^2 + \varepsilon(t) \left\| \frac{\partial^2 u^\varepsilon(x,t)}{\partial x^2} \right\|^2_{L^2(\Omega)} 
= \min\limits_{\begin{subarray}{c} s\in C^1(0,1) \end{subarray} } \frac{1}{n+1} \sum\limits^{n}\limits_{i=0} \left( s(x_i,t)-u^\delta_i \right)^2 + \varepsilon(t) \left\| \frac{\partial^2 s(x,t)}{\partial x^2} \right\|^2_{L^2(\Omega)} \\
\displaystyle   \leq \frac{1}{n+1} \sum\limits^{n}\limits_{i=0} \left( u(x_i,t)-u^\delta_i \right)^2 + \varepsilon(t) \left\| \frac{\partial^2 u(x,t)}{\partial x^2} \right\|^2_{L^2(\Omega)} \leq \delta^2 + \varepsilon(t) \left\| \frac{\partial^2 u(x,t)}{\partial x^2} \right\|^2_{L^2(\Omega)},
\end{array}
\end{align*}
which gives
\begin{equation}
\label{minimumIneq}
\left\| \frac{\partial u^\varepsilon(x,t)}{\partial x^2} \right\|^2_{L^2(\Omega)}  \leq \left\| \frac{\partial u(x,t)}{\partial x^2} \right\|^2_{L^2(\Omega)}.
\end{equation}
Consequently, we deduce, with the identity \eqref{Lemma1Eq}, that
\begin{equation} \label{equationbounded}
\left\Vert  \frac{\partial^2 e_1}{\partial x^2 }  \right\Vert \leq \left\Vert \frac{\partial^2 u^\varepsilon(x,t)}{\partial x^2 } \right\Vert + \left\Vert \frac{\partial^2 u(x,t)}{\partial x^2 } \right\Vert \leq 2 \left\Vert \frac{\partial^2 u(x,t)}{\partial x^2 } \right\Vert.
\end{equation}
Therefore, we obtain the following bound for $I_1$:
\begin{equation} \label{boundforI1}
I_1 \leq 2h\left\Vert \frac{\partial e_1}{\partial x } \right\Vert \, \left\Vert \frac{\partial^2 u(x,t)}{\partial x^2 }  \right\Vert.
\end{equation}

Next, we bound $I_2$ using the Cauchy–Schwarz inequality and \eqref{chiequation}. This yields
\begin{align}
I_2 ^2 \leq \sum_{i=1}^{n-1} e_1^2\left(x_i,t\right) \sum_{i=1}^{n-1}
\left(\chi_i-\chi_{i+1}\right)^2\nonumber  = \sum_{i=1}^{n-1}e_1^2\left(x_i,t\right)
\sum_{i=1}^{n-1}\frac{1}{h^2}\left(\int_{x_{i-1}}^{x_i}
\left(\frac{\partial e_1}{\partial x }\left(x\right)-\frac{\partial e_1}{\partial x }\left(x+h\right)\right)dx\right)^2. \nonumber
\end{align}
By construction,
\begin{align*}
\sum_{i=1}^{n-1} e_1^2\left(x_i,t\right) & = \sum_{i=1}^{n-1}
\left(u^\varepsilon \left(x_i,t\right)- u\left(x_i,t\right)\right)^2 \leq 2 \sum_{i=1}^{n-1} \left(
\left(u^\varepsilon \left(x_i,t\right)- u^\delta_i \right)^2 + \left(u^\delta_i- u\left(x_i,t\right)\right)^2 \right)
\\ & \leq 2 \sum_{i=0}^{n} \left(
u^\varepsilon \left(x_i,t\right)- u^\delta_i \right)^2 + 2(n-1) \delta^2 =
4n\delta^2
\end{align*}
and hence
\begin{align}
I_2 ^2 &\leq 4n\delta^2\sum_{i=1}^{n-1}\left(\int_{x_{i-1}}^{x_i}
\int_{x}^{x+h} \left \lvert \frac{\partial^2 e_1}{\partial x^2 } \left(\xi\right)\right \rvert d\xi \, dx \right)^2  / h^2 \nonumber
\\
&\leq
4n\delta^2
\sum_{i=1}^{n-1} \left(\int_{x_{i-1}}^{x_i} \int_{x_{i-1}}^{x_{i+1}}
\left \lvert \frac{\partial^2 e_1}{\partial x^2 }\left(\xi\right)\right \rvert d\xi \, dx \right)^2 / h^2
\leq
4n\delta^2\sum_{i=1}^{n-1}\left(\int_{x_{i-1}}^{x_{i+1}}\left \lvert \frac{\partial^2 e_1}{\partial x^2 } \left(\xi\right)
\right \rvert d\xi \right)^2 . \nonumber
\end{align}
This last integral can be bounded using Cauchy–Schwarz and \eqref{equationbounded} again (from the definition of $h$, $h>1/n$ holds):
\begin{align}
I_2^2 \leq
4n\delta^2\sum_{i=1}^{n-1}\int_{x_{i-1}}^{x_{i+1}}\left \lvert \frac{\partial^2 e_1}{\partial x^2 }\left(\xi\right)
\right \rvert^2d\xi \int_{x_{i-1}}^{x_{i+1}} d\xi \leq 16 \delta^2  \left\Vert \frac{\partial^2 e_1}{\partial x^2 }\right\Vert^2_{L^2(\Omega)} 
\leq 64\delta^2\left\Vert \frac{\partial^2 u(x,t)}{\partial x^2 } \right\Vert^2_{L^2(\Omega)}.
\nonumber
\end{align}

Inserting this and \eqref{boundforI1} into \eqref{normforderivativeofe}, we finally obtain
\begin{equation} \label{estimfornormofe}
\left\Vert \frac{\partial e_1}{\partial x } \right\Vert^2_{L^2(\Omega)} \leq 2h\left\Vert \frac{\partial e_1}{\partial x } \right\Vert_{L^2(\Omega)} \, \left\Vert \frac{\partial^2 u(x,t)}{\partial x^2 }  \right\Vert_{L^2(\Omega)} +8\delta\left\Vert \frac{\partial^2 u(x,t)}{\partial x^2 }
\right\Vert_{L^2(\Omega)}.
\end{equation}
Completing the squares permits us to conclude from \eqref{estimfornormofe} that
\begin{align}
\left(\left\Vert \frac{\partial e_1}{\partial x } \right\Vert_{L^2(\Omega)} -h \left\Vert \frac{\partial^2 u(x,t)}{\partial x^2 }  \right\Vert_{L^2(\Omega)} \right)^2 
\leq \left(h\left\Vert \frac{\partial^2 u(x,t)}{\partial x^2 }  \right\Vert_{L^2(\Omega)} +
\sqrt{8}\sqrt{\delta}\left\Vert \frac{\partial^2 u(x,t)}{\partial x^2 }  \right\Vert^{1/2}_{L^2(\Omega)} \right)^2.
\nonumber
\end{align}
This yields
\begin{equation}
\left\Vert \frac{\partial e_1}{\partial x } \right\Vert_{L^2(\Omega)} \leq 2h \left\Vert \frac{\partial^2 u(x,t)}{\partial x^2 }  \right\Vert_{L^2(\Omega)} + \sqrt{8}\sqrt{\delta}\left\Vert \frac{\partial^2 u(x,t)}{\partial x^2 }
\right\Vert^{1/2}_{L^2(\Omega)},
\nonumber
\end{equation}
and inequalities \eqref{PoincareIneq} and \eqref{Lemma2Ineq} imply that, for a.e. $t\in \bar{\mathcal{T}}$,
\begin{align*}
\begin{split}
& \Vert e(\cdot,t) \Vert_{H^1(0,1)} \leq 5\left\Vert \partial_x e(x,t) \right\Vert_{L^2(\Omega)} \leq 5 \left( \left\Vert \frac{\partial e_1}{\partial x } \right\Vert_{L^2(\Omega)} + \left\Vert \frac{\partial e_2}{\partial x } \right\Vert_{L^2(\Omega)} \right)  \\ & \quad
\leq 5\left(2h \left\Vert \frac{\partial^2 u(x,t)}{\partial x^2 }  \right\Vert_{L^2(\Omega)} + \sqrt{8}\sqrt{\delta}\left\Vert \frac{\partial^2 u(x,t)}{\partial x^2 }
\right\Vert^{1/2}_{L^2(\Omega)} + \frac{h}{\pi} \left\|\frac{\partial^2 u(x,t)}{\partial x^2 }\right\|_{L^2(\Omega)} \right),
\end{split}
\end{align*}
which yields the required estimate \eqref{NoisyErrIneq}.
\end{proof}

\section{Numerical examples}
\label{simulation}

In this section, we present some numerical experiments to illustrate the efficiency of our new approach. For each example, we first verify the numerical behavior of the asymptotic solution, whose accuracy is theoretically guaranteed by Theorem \ref{MainThm}, and then subsequently demonstrate the efficiency of Algorithm 1 for the corresponding inverse problems.

We consider the following reaction–diffusion–advection equation with source function $f(x)$, which will be set differently in the simulation study:
\begin{align} \label{forwardexample1}
\begin{cases}
\displaystyle \mu \frac{\partial^2 u}{\partial x^2} - \frac{\partial u}{\partial t} = -u \frac{\partial u}{\partial x}+f(x), \quad  x \in (0, 1), \quad  t \in (0,T], \quad \mu=0.01, \\
\displaystyle  u(0,t)=u^l, \quad u(1,t)=u^r, \quad t \in [0,T], \\
\displaystyle  u(x,0)=u_{init}, \quad x \in [0, 1].
\end{cases}
\end{align}
According to Theorem \ref{MainThm}, we need to verify Assumptions \ref{A1}–\ref{A4}, which stand for the sufficiency conditions for the existence of an asymptotic solution to problem \eqref{forwardexample1}. We therefore repeat the procedure presented in Subsection \ref{constructionOfAsymptotic}. By solving two equations \eqref{LeftRightSolution} we obtain the main regular terms $\varphi^l (x)$ and $ \varphi^r (x)$ in the forward problems.
The problem for determining the leading term of the asymptotic description of the front $x_0 (t)$ takes the form
\begin{align}\label{eq44Ex1}
\begin{cases}
\displaystyle \frac{dx_0(t)}{dt}=- \frac{1}{2}\left( \varphi^l (x_0(t))+ \varphi^r (x_0(t)) \right),\\
x_0(0)=0.1.
\end{cases}
\end{align}

Then, by solving \eqref{eq44Ex1} numerically, we can verify Assumption \ref{A3}.

We take the initial function in the form
$$\displaystyle u_{init} (x,\mu )=\frac{u^r-u^l}{2} \tanh\left(\frac{x-0.1}{\mu}\right)+ \frac{u^r+u^l}{2},$$
with an inner transition layer in the vicinity of $x=0.1$.

 Thus, if Assumptions \ref{A1}–\ref{A4} are satisfied, the considered equation \eqref{forwardexample1} has the following solution:
\begin{align} \label{asymptoticsolEXAMPLE1}
U_{0}(x,t)=\begin{cases}
\displaystyle \varphi^l (x) +\frac{ \left(\varphi^{r}(x_{0}(t))-\varphi^{l}(x_{0}(t)) \right)}{  \exp \left( \left(x-x_0(t) \right) \left( \frac{\varphi^{l}(x_{0}(t))-\varphi^{r}(x_{0}(t)) }{2\mu} \right) \right)+ 1}  ,  \ x \in [0;x_0(t)],\\
\displaystyle \varphi^r (x)+\frac{ -\left(\varphi^{r}(x_{0}(t))-\varphi^{l}(x_{0}(t))\right)}{  \exp \left(\left(x-x_0(t) \right)  \left( \frac{\varphi^{r}(x_{0}(t))-\varphi^{l}(x_{0}(t)) }{2\mu}\right) \right)+ 1} ,\  x \in [x_0(t);1].
\end{cases}
\end{align}

For the simulation of inverse problems \textbf{(IP)}, we consider the problem of identifying the source function $f(x)$ in the nonlinear PDE model \eqref{forwardexample1}. The numerical experiments consist of three steps. First, we obtain the synthetic exact measurement data $\{u(x_i,t_0), \frac{\partial u}{\partial x}(x_i,t_0)\}^n_{i=0}$ by solving  the forward problem  \eqref{forwardexample1} numerically with the finite volume method, where we introduce a mesh uniformly with respect to spatial variable $\Theta=\left\{ x_i, 0 \leq i \leq n: x_i= h i, h = 1/n \right\} $. Second, we generate the artificial noisy data by  adding independent and identically distributed (i.i.d.) random variables with a uniform distribution with noise level $\delta$; i.e. for $i=0,\cdots, n$,
\begin{equation}
\label{noiseUni}
u^{\delta}_i := [1+\delta(2\, \text{rand} -1)] u(x_i,t_0), ~ w^{\delta}_i := [1+\delta(2\, \text{rand} -1)] \frac{\partial u}{\partial x}(x_i,t_0),
\end{equation}
where $\text{rand}$ returns a pseudo-random value drawn from a uniform distribution on $[0,1]$.
In the last simulation step, the observed data is processed by Algorithm \ref{alg:Framwork}, and the retrieved source function is compared with the one from the input. Moreover, we also output the relative a posteriori errors of the estimated source function and the lower and upper source functions.

\subsection{Example 1}
\subsubsection{Forward problem}
\begin{figure}[H]
\begin{center}
\includegraphics[width=0.4\linewidth]{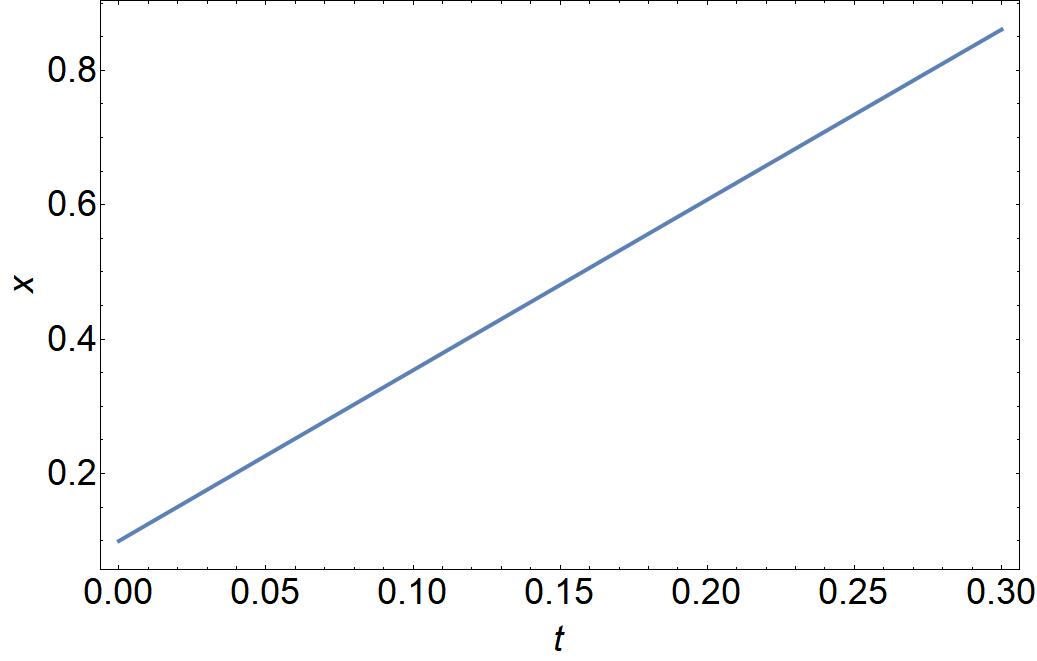}
\caption{Numerical solution of \eqref{eq44Ex1} for $t \in [0,  0.3] $.}
\label{fig:inverseMeshEX1}
\end{center}
\end{figure}

\begin{figure}[H]
\vspace{-3ex} \centering
\subfigure[]{
\includegraphics[width=0.4\linewidth]{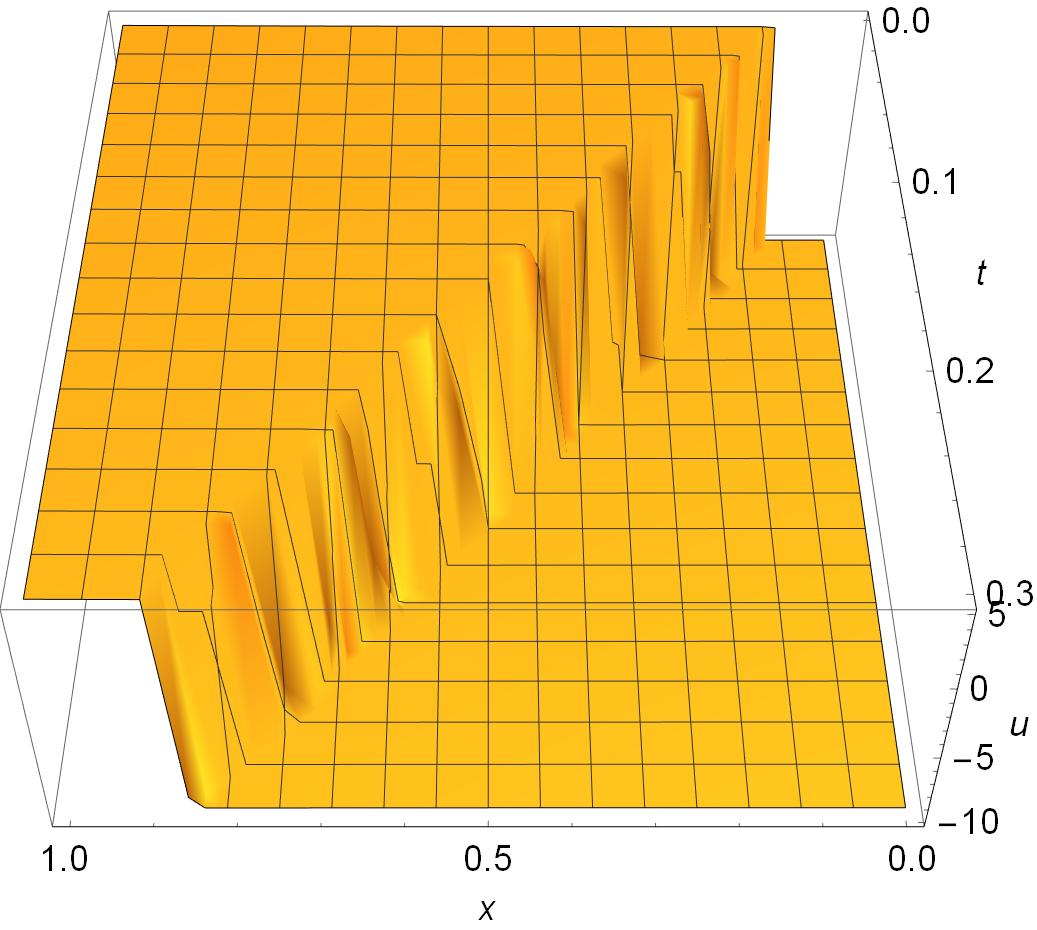} \label{fig:asymptoticsolexample1EX1} }
\hspace{0ex}
\subfigure[]{
\includegraphics[width=0.4\linewidth]{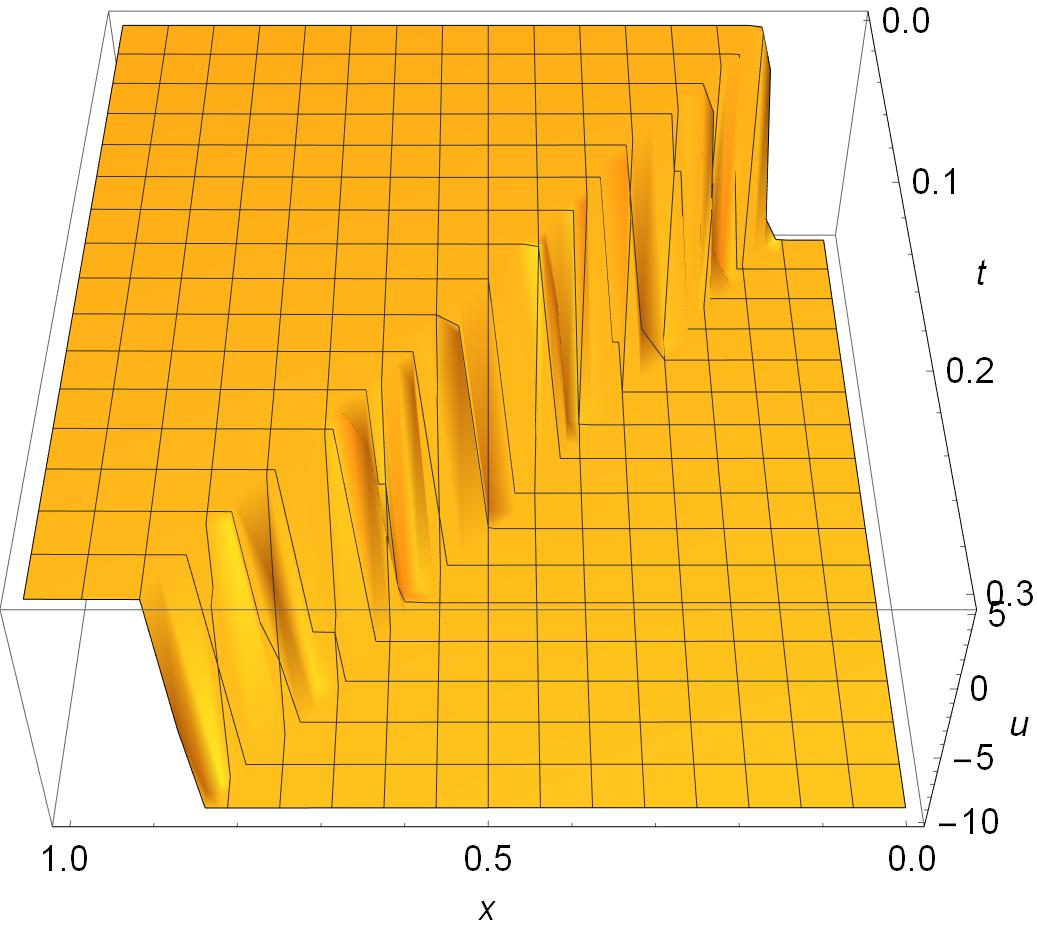} \label{fig:numericsolexample1EX1} }
\caption{Asymptotic solution \subref{fig:asymptoticsolexample1EX1} and numerical solution (using the finite-volume method) \subref{fig:numericsolexample1EX1} of PDE \eqref{forwardexample1}  with $f^*(x)=x-x^2+x^3$, $x\in [0,1]$, $t \in [0,0.3],\mu =0.01$.} \label{solutionsexample1EX1}
\end{figure}

In this example, we consider equation \eqref{forwardexample1} with a given monotonically increasing source function $f^*(x)=x-x^2+x^3$ and parameters $u^l=-10, u^r=5, T=0.3$. We explicitly find the zero-order regular functions,
\begin{align*}
\varphi^l (x)= - \frac{\sqrt{600+6x^2-4x^3+3x^4}}{ \sqrt{6}}, \quad
\varphi^r (x) = \frac{\sqrt{145+6x^2-4x^3+3x^4}}{ \sqrt{6}},
\end{align*}
and numerically verify that $0<x_0(t)<1$ for all $t\in \bar{\mathcal{T}}$ (Fig. \ref{fig:inverseMeshEX1}). The initial function takes the form $\displaystyle u_{init} (x,\mu )=7.5\tanh\left(\frac{x-0.1}{0.01}\right)-2.5$. Thus, Assumptions \ref{A1}–\ref{A4} are satisfied, and the asymptotic solution is shown in Fig. \ref{fig:asymptoticsolexample1EX1}. We also draw the numerical solution (using the finite-volume method) for problem \eqref{forwardexample1} in Fig. \ref{fig:numericsolexample1EX1}, which will be used as the high resolution of the exact solution $u$. The relative error of the asymptotic solution is $ \frac{\| U_0(x,t) - u(x,t) \|_{L^{2}(\bar{\Omega} \times [0,0.3])}}{\|u(x,t) \|_{L^{2}(\bar{\Omega} \times [0,0.3])}} = 0.0586$.

\subsubsection{Inverse problem} \label{Ex1InversePr}

In the simulation, we use the error level $\delta = 1 \%$, $t_0=0.2$, $n^l=12$, $n^r=13$, and $n=20$, and take the values of the grid $u(x_i,t_0)$ and  $\frac{ \partial u}{\partial x} (x_i,t_0)$ from the forward problem.

We skip the points from transition layer $(x_0(t_0)-\Delta x /2, x_0(t_0)+\Delta x/2)$, and use nodes in only two regions, located on two sides of the transition layer, with indices $i= 0, \cdots, n^l$ and $i=n^r, \cdots, n$. The uniform noise \eqref{noiseUni} is added to the values $u(x_i, t_0)$ and  $\frac{ \partial u }{\partial x} (x_i, t_0)$ to produce noisy data  $\{u^{\delta}_i, w^{\delta}_i\}^{n^l}_{i=0}$ and $\{u^{\delta}_i, w^{\delta}_i\}^{n}_{i=n^r}$ on the left and right intervals with respect to the transition layer.

Following Algorithm \ref{alg:Framwork},  we
obtain the approximate source function by solving the following optimization problem:
\begin{equation}
\label{uAlphaEx1}
f^\delta(x) = \mathop{\arg\min}_{\begin{subarray}{c} f\in C^1(0,1): \\   f(x_0)<f(x_1)<\cdots<f(x_n)  \end{subarray} } \frac{1}{n+1} \sum^{n}_{i=0} \left( f(x_i)-u^{\delta}_i w^{\delta}_i \right)^2.
\end{equation}

\begin{figure}[H]
\begin{center}
\includegraphics[width=0.4\linewidth]{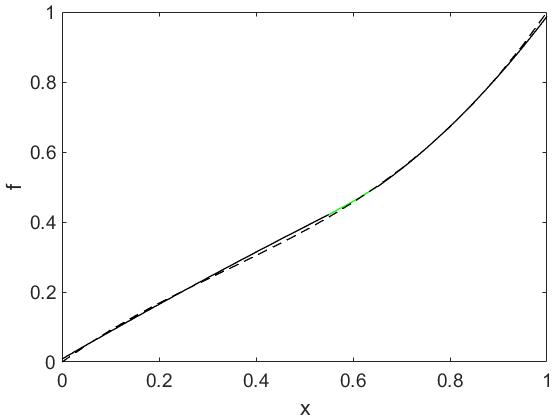}
\caption{The result of reconstructing the source function $f^\delta(x)$ (black lines) for $t_0 = 0.2$. The green line interpolates the source function in the transition layer; this can be compared with the exact source function $f^*(x)=x-x^2+x^3$ (dashed line). }
\label{fig:sourseFuncEx1}
\end{center}
\end{figure}

The reconstructed source function $f^\delta(x)$ is shown in Fig. \ref{fig:sourseFuncEx1}.

The relative error of the reconstruction is $\|f^\delta - f^* \|_{L^2(0,1)} / \|  f^* \|_{L^2(0,1)} = 0.0076$.

Using formula \eqref{posterioriErr2}, we can calculate the relative a posteriori error for the obtained approximate source function $\Delta_1=0.0641$, which is slightly larger than the value of the relative error.

\begin{figure}[H]
\begin{center}
\includegraphics[width=0.4\linewidth]{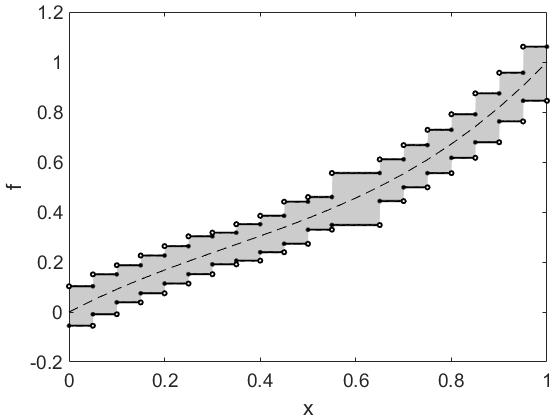}
\caption{ The lower $f^{low}(x)$ and upper $f^{up}(x)$ solutions, which may be compared with the accuracy source function $f^*(x)=x-x^2+x^3$ (dashed line). }
\label{fig:monotonicUPLOW}
\end{center}
\end{figure}

Formulas \eqref{monotonicLOW} and \eqref{monotonicUP}, for monotonic functions, were used to construct the lower $f^{low}(x)$ and upper $f^{up}(x)$ solutions, whose figures are shown in Fig. \ref{fig:monotonicUPLOW}. In this figure, we can see that the exact source function $f^*$ lies between these two functions, i.e. the underground truth $f^*$ is located in the shadow region in Fig. \ref{fig:monotonicUPLOW}.

\subsection{Example 2}
\subsubsection{Forward problem}
\begin{figure}[H]
\begin{center}
\includegraphics[width=0.4\linewidth]{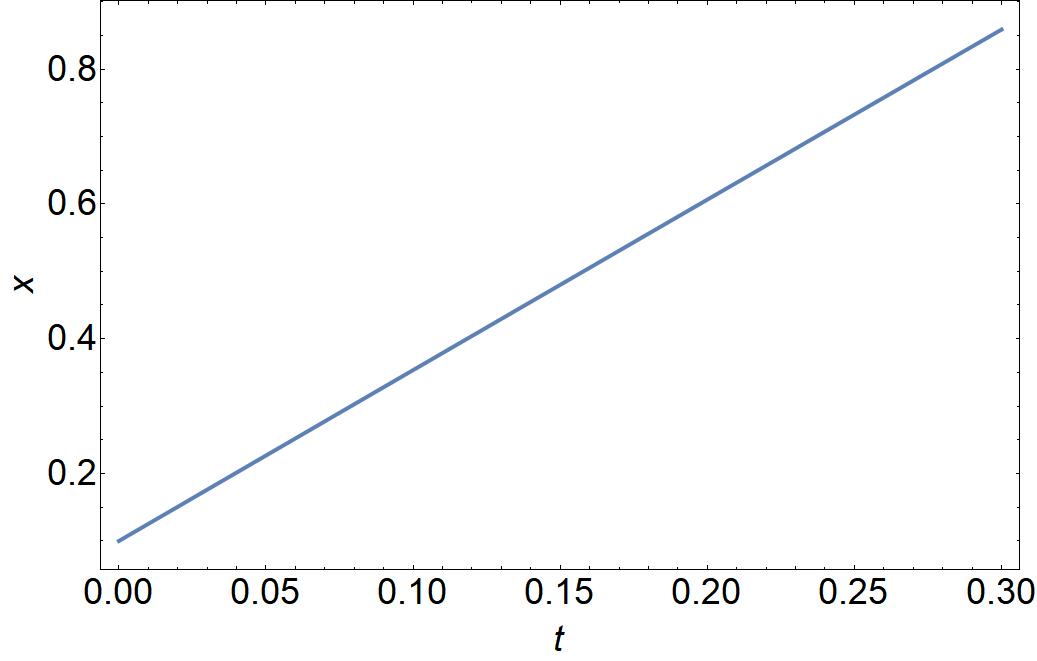}
\caption{Numerical solution of \eqref{eq44Ex1} for $f^*(x)=\sqrt{x-x^2} $, $t \in [0,  0.3] $.}
\label{fig:inverseMeshEX2}
\end{center}
\end{figure}

In this example, we consider PDE \eqref{forwardexample1} with a given convex source function $f^*(x)=\sqrt{x-x^2}$ and parameters $u^l=-10, u^r=5, T=0.3$.
The regular functions of zero order have the form
$$\varphi^l(x)=-\frac{1}{2} \sqrt{(4 x-2) \sqrt{(1-x) x}-2 \sin ^{-1}\left(\sqrt{1-x}\right)+\pi +400},$$
$$\varphi^r (x)= \frac{1}{\sqrt{2}} \sqrt{(2 x-1) \sqrt{(1-x) x}-\sin ^{-1}\left(\sqrt{1-x}\right)+50},$$
and we numerically verify that $0<x_0(t)<1$ for all $t\in[0,0.3]$ (see Fig. \ref{fig:inverseMeshEX2}).

\begin{figure}[H]
\vspace{-3ex} \centering
\subfigure[]{
\includegraphics[width=0.4\linewidth]{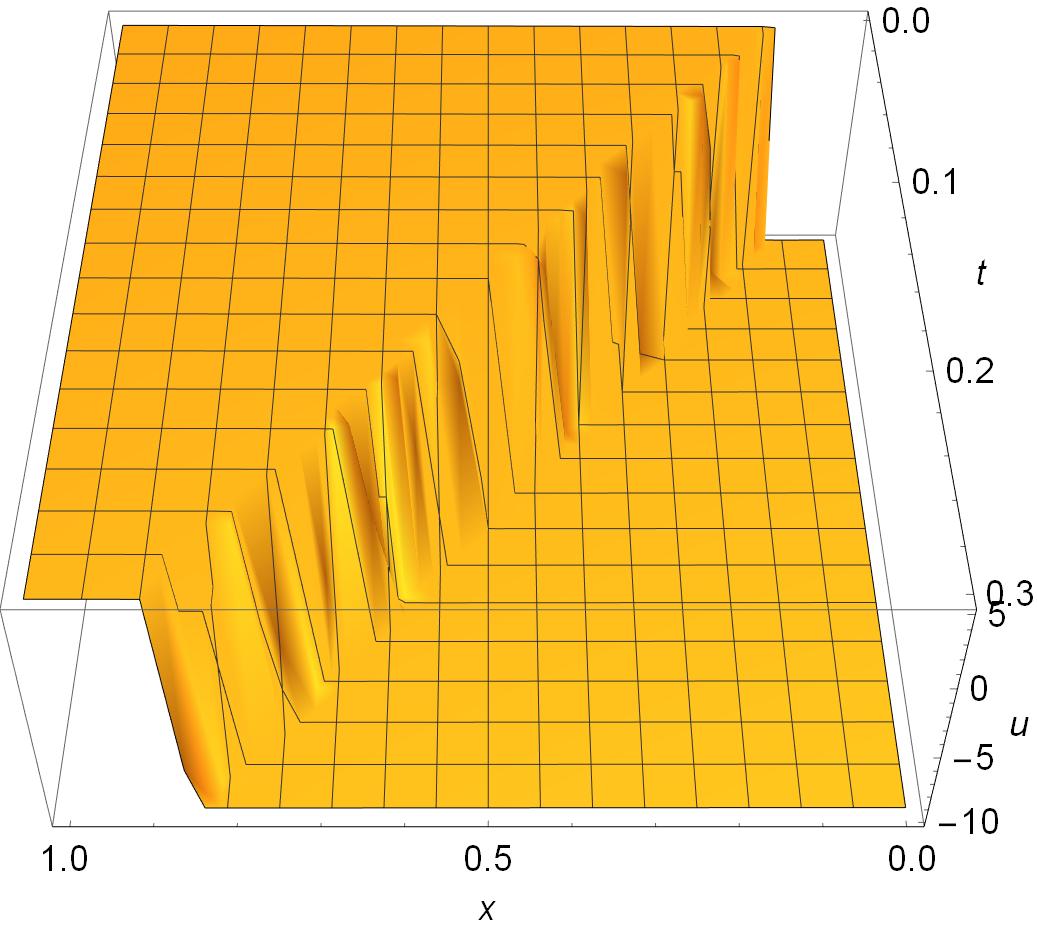} \label{fig:asymptoticsolexample1EX2} }
\hspace{0ex}
\subfigure[]{
\includegraphics[width=0.4\linewidth]{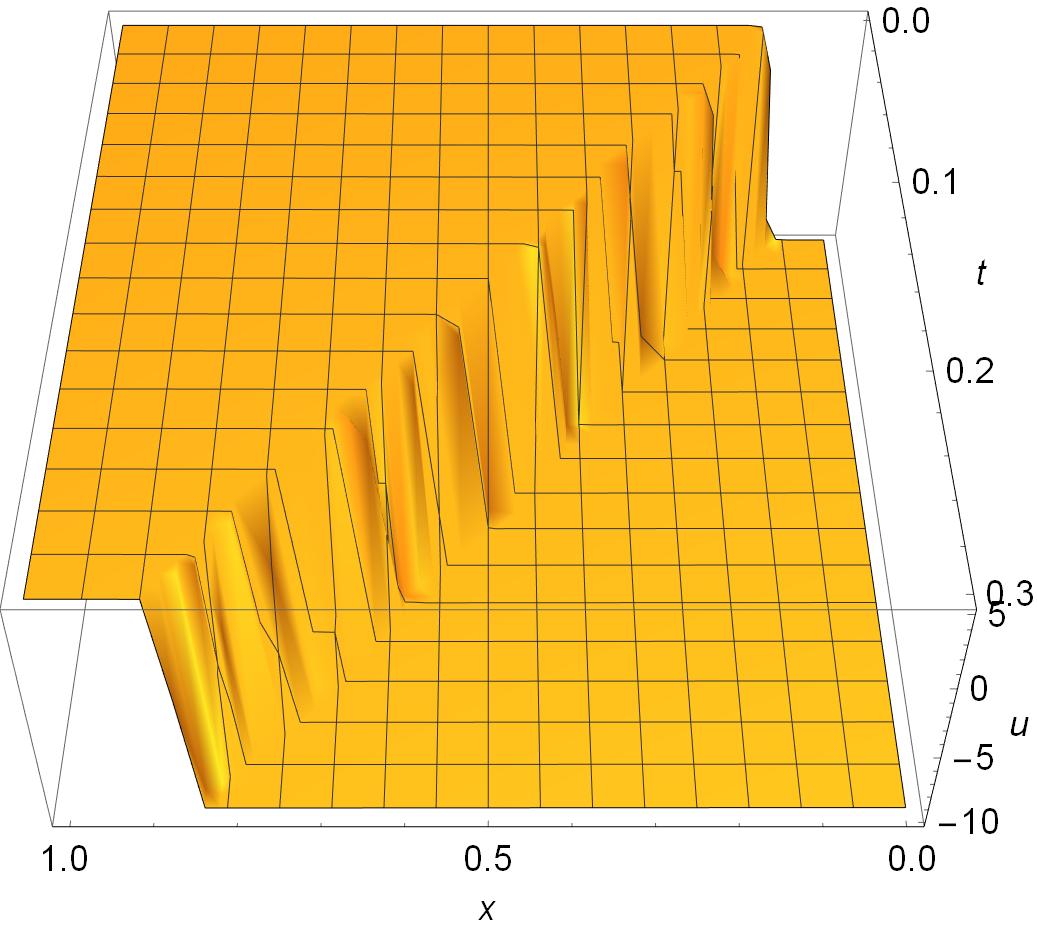} \label{fig:numericsolexample1EX2} }
\caption{Asymptotic solution \subref{fig:asymptoticsolexample1EX1} and numerical solution (using the finite-volume method) \subref{fig:numericsolexample1EX2} of PDE \eqref{forwardexample1}  with $f^*(x)=\sqrt{x-x^2}$, $x\in [0,1]$, $t \in [0,0.3],\mu =0.01$.} \label{solutionsexample1EX2}
\end{figure}

The initial function takes the form $\displaystyle u_{init} (x,\mu )=7.5\tanh\left(\frac{x-0.1}{0.01}\right)-2.5$. Thus, Assumptions \ref{A1}–\ref{A4} are fulfilled, and the asymptotic solution is shown in Fig. \ref{fig:asymptoticsolexample1EX2}. We also draw the numerical solution (using the finite-volume method) for problem \eqref{forwardexample1} in Fig. \ref{fig:numericsolexample1EX2}. The relative error of the asymptotic solution is $ \frac{\| U_0(x,t) - u(x,t) \|_{L^{2}(\bar{\Omega} \times [0,0.3])}}{\|u(x,t) \|_{L^{2}(\bar{\Omega} \times [0,0.3])}} = 0.0386$.

\subsubsection{Inverse problem}

\begin{figure}[H]
\begin{center}
\includegraphics[width=0.4\linewidth]{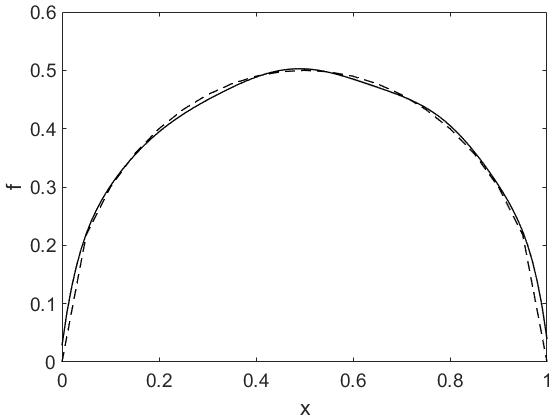}
\caption{The result of reconstructing the source function $f^\delta(x)$ (solid black line) for $t_0 = 0.2 $; this can be compared with the exact source function $f^*(x)=\sqrt{x-x^2}$ (dashed line). }
 \label{fig:sourseFuncEx2}
\end{center}
\end{figure}

For the inverse problem, we use the same input settings as in the inverse problem from Example 1. According to our theoretical analysis (e.g. Theorem 1), we can exclude data values belonging to the transition layer, and then produce noisy data  $\{u^{\delta}_i\}^{n^l}_{i=0}, \{w^{\delta}_i\}^{n^l}_{i=0}$ and  $\{u^{\delta}_i\}^{n}_{i=n^r}, \{w^{\delta}_i\}^{n}_{i=n^r}$. The approximate source function is estimated by solving the optimization problem
\begin{equation}
\label{uAlphaEx2}
f^\delta(x) = \mathop{\arg\min}_{\begin{subarray}{c} f\in C^1(0,1):\\  f(x_{i-1})-2f(x_i)+f(x_{i+1}) >0, \\ 1 \leq i < n-1\end{subarray} } \frac{1}{n+1} \sum^{n}_{i=0} \left( f(x_i)-u^{\delta}_i w^{\delta}_i  \right)^2,
\end{equation}
and the optimization result is shown in Fig. \ref{fig:sourseFuncEx2}.

The relative error of the recovered source function is $\|f^\delta - f^* \|_{L^2(0,1)} / \|  f^* \|_{L^2(0,1)} = 0.0297$.

Using formula \eqref{posterioriErr2}, we calculate the relative a posteriori error for the obtained approximate source function $\Delta_1=0.2552$.

\begin{figure}[H]
\begin{center}
\includegraphics[width=0.4\linewidth]{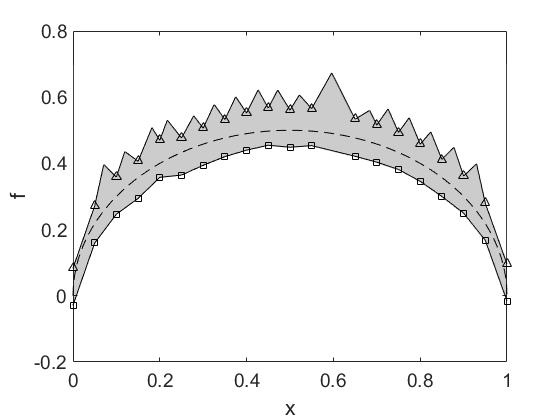}
\caption{The lower $f^{low}(x)$ (line with squares) and upper $f^{up}(x)$ (line with triangles) solutions; these can be compared with the exact source function  $f^*(x)=\sqrt{x-x^2} $ (dashed line). }
\label{fig:convexUPLOW}
\end{center}
\end{figure}

Fig. \ref{fig:convexUPLOW} shows the lower $f^{low}(x)$ and upper $f^{up}(x)$ solutions, which are constructed according to formulas \eqref{convexLOW} and \eqref{convexUP} for convex functions. Fig. \ref{fig:convexUPLOW} also indicates that the exact source function is located in the shadowy area between the upper and lower solutions.

\subsection{Example 3}
\subsubsection{Forward problem}
In the last example, we consider equation \eqref{forwardexample1} with a given source function $f^*(x)=x\sin(3 \pi x)$ and parameters $u^l=-8, u^r=4, T=0.2$.

\begin{figure}[H]
\begin{center}
\includegraphics[width=0.4\linewidth,height=0.6\textwidth,keepaspectratio]{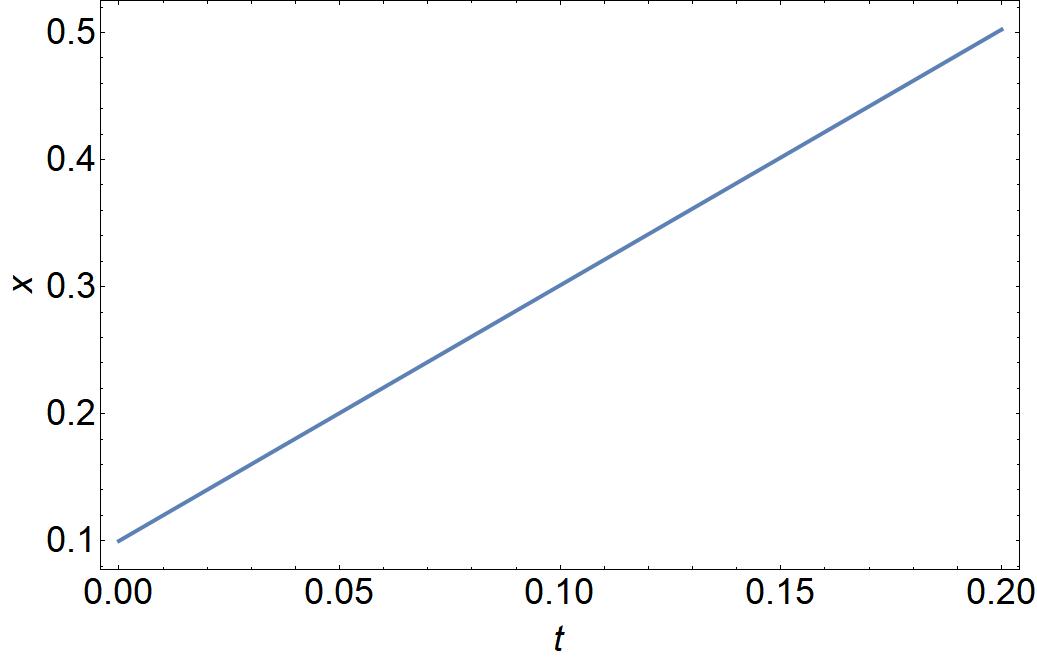}
\caption{Numerical solution of \eqref{eq44Ex1} for $f^*(x)=x\sin(3 \pi x)$, $t \in [0,  0.2] $.}
\label{fig:x0example1}
\end{center}
\end{figure}

We explicitly find the zero-order regular functions,
\begin{align*}
\varphi^l (x)= - \frac{\sqrt{2(288 \pi^2 - 3 \pi x\cos(3 \pi x)+\sin(3 \pi x))}}{3 \pi },\\
\varphi^r (x) = \frac{\sqrt{2(72 \pi^2 -3 \pi - 3 \pi x\cos(3 \pi x)+\sin(3 \pi x))}}{3 \pi },
\end{align*}
and numerically verify that $0<x_0(t)<1$ for all $t\in [0,0.2]$ (see Fig. \ref{fig:x0example1}).

The initial function has the form $\displaystyle u_{init} (x,\mu ) = 6\tanh \left(\frac{x-0.1}{0.01}\right) -2$. Thus, Assumptions \ref{A1}–\ref{A4} are verified and the considered equation \eqref{forwardexample1} for a given source function has asymptotic solution in the form of an autowave with a transitional moving layer localized near  $x_0(t)$, which for $t \in [0, 0.2]$ in the zero approximation has the form shown below in Fig. \ref{fig:asymptoticsolexample1}. The numerical solution of \eqref{forwardexample1} using the finite volume is displayed in Fig. \ref{fig:numericsolexample1}.

\begin{figure}[H]
\vspace{-3ex} \centering
\subfigure[]{
\includegraphics[width=0.4\linewidth]{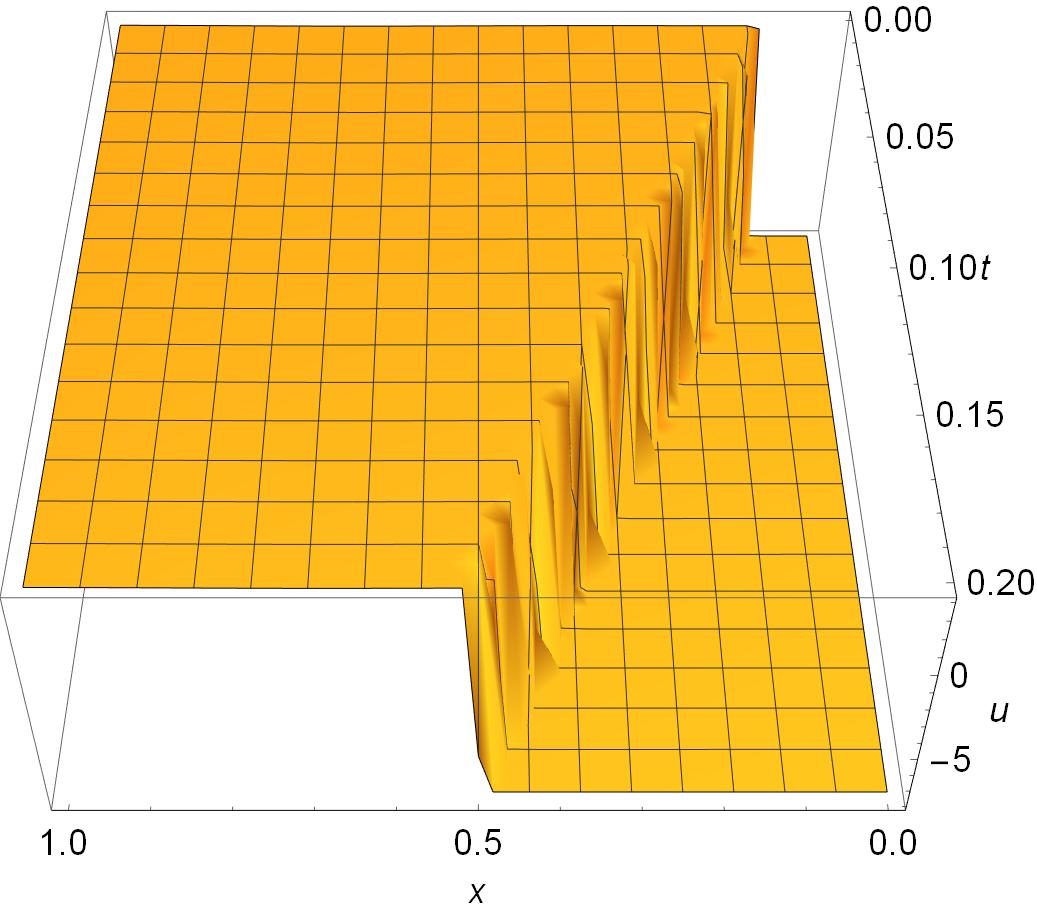} \label{fig:asymptoticsolexample1} }
\hspace{0ex}
\subfigure[]{
\includegraphics[width=0.4\linewidth]{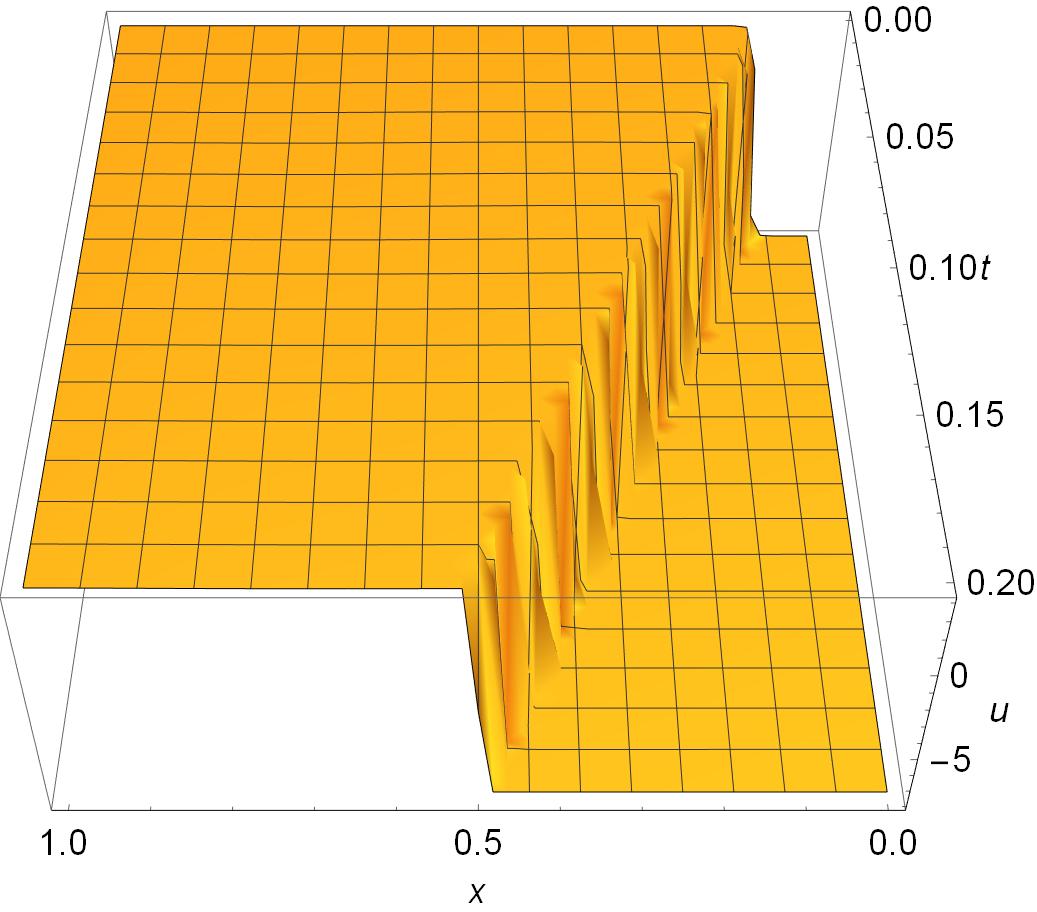} \label{fig:numericsolexample1} }
\caption{Asymptotic solution \subref{fig:asymptoticsolexample1} and numerical solution (using the finite-volume method) \subref{fig:numericsolexample1} of PDE \eqref{forwardexample1}  with $f(x)=x\sin(3 \pi x)$, $x\in [0,1]$, $t \in [0,0.2], \mu =0.01$.} \label{solutionsexample1}
\end{figure}

The relative error of the asymptotic solution is
$$ \frac{\| U_0(x,t) - u(x,t) \|_{L^{2}(\bar{\Omega} \times [0,0.2])}}{\|u(x,t) \|_{L^{2}(\bar{\Omega} \times [0,0.2])}} = 0.0411. $$

Note that, even if we use only the regular part of asymptotics \eqref{u0regu} as a solution to problem \eqref{forwardexample1}, because of the narrow transition layer, the relative error is still small:
$$ \frac{\| \bar{u}_{0}(x) - u(x,t) \|_{L^{2}(\bar{\Omega} \times [0,0.2])}}{\|u(x,t) \|_{L^{2}(\bar{\Omega} \times [0,0.2])}} = 0.1081. $$

\subsubsection{The inverse source problem}
\begin{figure}[H]
\begin{center}
\includegraphics[width=0.3\linewidth]{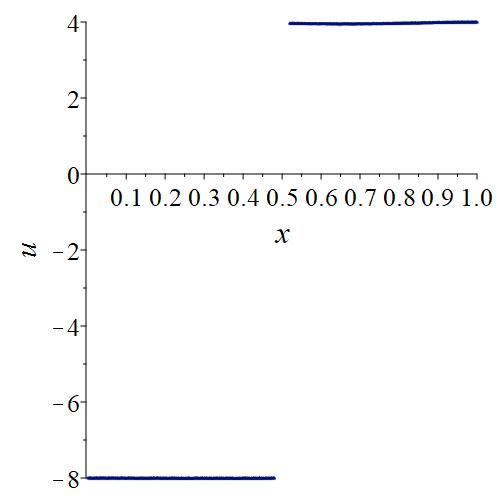}
\caption{ Solution values $u(x_i, t_0)$ on the left $i= 0, \cdots, n^l$ and on the right $i=n^r, \cdots, n$ intervals for $t_0 = 0.2$.}
\label{fig:inverseMesh}
\end{center}
\end{figure}

Now we consider the problem of identifying the source function $f(x)$ in the previously described PDE model \eqref{forwardexample1}. For this example, we assume that we know only the values of $\{u(x_i,t_0)\}^n_{i=0}$ at time $t_0$. As in the previous examples, synthetic measurement data are obtained from the numerical result using the finite-volume method for the forward problem \eqref{forwardexample1} (see Fig. \ref{fig:numericsolexample1}). We introduce a mesh uniformly with respect to spatial variable $\Theta=\left\{ x_i, 0 \leq i \leq n: x_i= h i, h = 1/n \right\}$, and use nodes in only two regions outside the transition layer, i.e. $[0,x_0(t_0)-\Delta x /2 ]$ and $[x_0(t_0)+\Delta x/2 ,1]$ with node indices $i= 0, \cdots, n^l$ and $i=n^r, \cdots, n$ (see Fig. \ref{fig:inverseMesh}). The i.i.d. uniform noises \eqref{noiseUni} with two noise levels $\delta_{i}$ ($i=1,2$) are added to $\{u(x_i,t_0)\}^n_{i=0}$ to produce noisy data  $\{u^{\delta_{1,2}}_i(t_0)\}^{n^l}_{i=0}$ and $\{u^{\delta_{1,2}}_i(t_0)\}^n_{i=n^r}$ on the left and right intervals, respectively.

In the simulation, we use parameters $\delta_1 = 0.1 \%$,  $\delta_2 = 1\%$, $t_0=0.2$, $n^l=238$, $n^r=260 $, and $n=499$. Following Algorithm \ref{alg:Framwork}, we obtain the smoothed data $u^\varepsilon(x, t_0)$ according to the following optimization problem for the left and right segments, respectively:
\begin{equation*}
\label{uAlphaEx3Left}
u^\varepsilon(x,t_0) = \mathop{\arg\min}_{\begin{subarray}{c} s\in C^1(0, h n^l) \end{subarray} } \frac{1}{n^l+1} \sum^{n^l}_{i=0} \left( s(x_i,t_0)-u^\delta_i \right)^2 + \varepsilon^l(t_0) \left\| \frac{\partial^2 s(x,t_0)}{\partial x^2} \right\|^2_{L^2(0, h n^l)},
\end{equation*}
\begin{equation*}
\label{uAlphaEx3Right}
u^\varepsilon(x,t_0) = \mathop{\arg\min}_{\begin{subarray}{c} s\in C^1(h n^r, 1) \end{subarray} } \frac{1}{n+1-n^r} \sum^{n}_{i=n^r} \left( s(x_i,t_0)-u^\delta_i \right)^2 + \varepsilon^r(t_0) \left\| \frac{\partial^2 s(x,t_0)}{\partial x^2} \right\|^2_{L^2(h n^r, 1)}.
\end{equation*}

\begin{figure}[H]
\vspace{-1ex} \centering
\subfigure[]{
\includegraphics[width=0.4\linewidth]{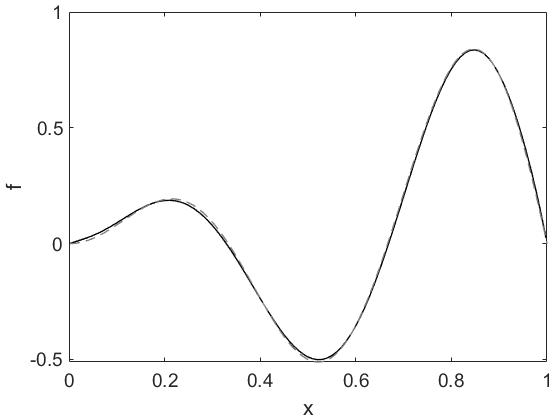} \label{fig:sourcedelta10-4} }
\hspace{0ex}
\subfigure[]{
\includegraphics[width=0.4\linewidth]{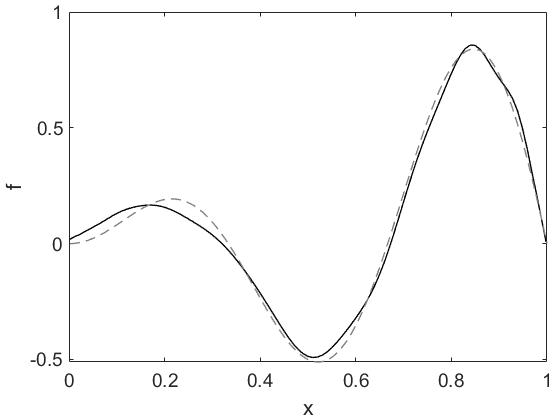} \label{fig:sourcedelta10-3} }
\caption{The results of reconstructing the source function $f(x)$ (black lines) for different input error levels $\delta_1=0.1 \%$ \subref{fig:sourcedelta10-4} and $\delta_2=1 \%$ \subref{fig:sourcedelta10-3}; they can be compared with the exact source function $f^*(x)=x \sin{(3 \pi x)}$ (dashed lines).} \label{fig:sourcereconstruction}
\end{figure}

According to the second step in Algorithm \ref{alg:Framwork}, we obtain the missing measurements of $\{w^{\delta}_i\}^{n^l}_{i=1} \cup \{w^{\delta}_i\}^{n}_{i=n^r+1}$ by taking the numerical derivative, i.e. $w^{\delta}_i = \frac{ u^\varepsilon(x_i,t_0)-u^\varepsilon(x_{i-1},t_0)}{x_i-x_{i-1}}$, for $1<i \leq n^l$ and $n^r+1<i \leq n$. Then, the regularized approximate source function $f^\delta(x)$ is computed by formula \eqref{uAlpha1}. The results are shown in Fig. \ref{fig:sourcereconstruction}, and from them we can conclude that our approach is stable and accurate.

The relative error and the relative a posteriori error of the reconstructed source functions $f^\delta$ for two different sets of noisy data are as follows:
\begin{itemize}
\item For $\delta_1 =0.1\%$: $\| f^\delta- f^* \|_{L^2(0,1)} / \|  f^* \|_{L^2(0,1)} = 0.0233$ and $\Delta_1=5.7082$;
\item For $\delta_2 =1\%$: $\| f^\delta- f^* \|_{L^2(0,1)} / \|  f^* \|_{L^2(0,1)}= 0.0894$ and $\Delta_1=6.4588$.
\end{itemize}
They indicate that, for this model problem, our relative a posteriori error $\Delta_1$ is slightly over-estimated. Nevertheless, the reasonable value of $\Delta_1$ is always useful in practice for real-world problems.

\begin{figure}[H]
\begin{center}
\includegraphics[width=0.4\linewidth]{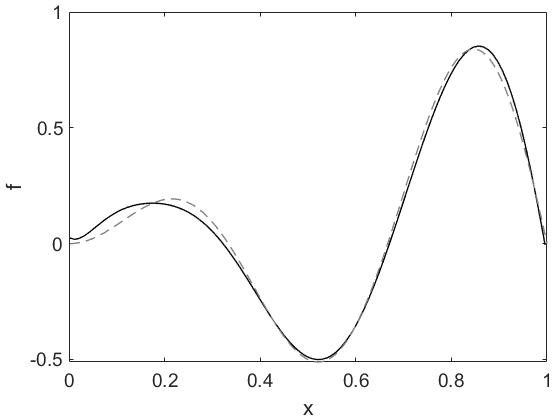}
\caption{ The result of reconstructing the source function $f(x)$ (solid black line) for the input-data error level $\delta_2=1 \%$; this can be compared with the exact source function $f^*(x)=x \sin{(3 \pi x)}$ (dashed line).}
\label{fig:sourcedeltaGAP10-3}
\end{center}
\end{figure}

Note that even if the initial data has missing points we still can reconstruct the source function. In Fig. \ref{fig:sourcedeltaGAP10-3} we reconstruct the source function for $t_0=0.17$ when the initial data (with 1\% noise level) on the right interval has a gap between the points $x=0.77$ and $x=0.87$, and missing values for the noisy data were approximated by a first-degree spline. In this case, the
relative error of the reconstruction equals $\|f^\delta - f^* \|_{L^2(0,1)} / \|  f^* \|_{L^2(0,1)} = 0.0722 $, while the relative a posteriori error is $\Delta_1= 6.25$.

\begin{figure}[H]
\begin{center}
\includegraphics[width=0.5\linewidth]{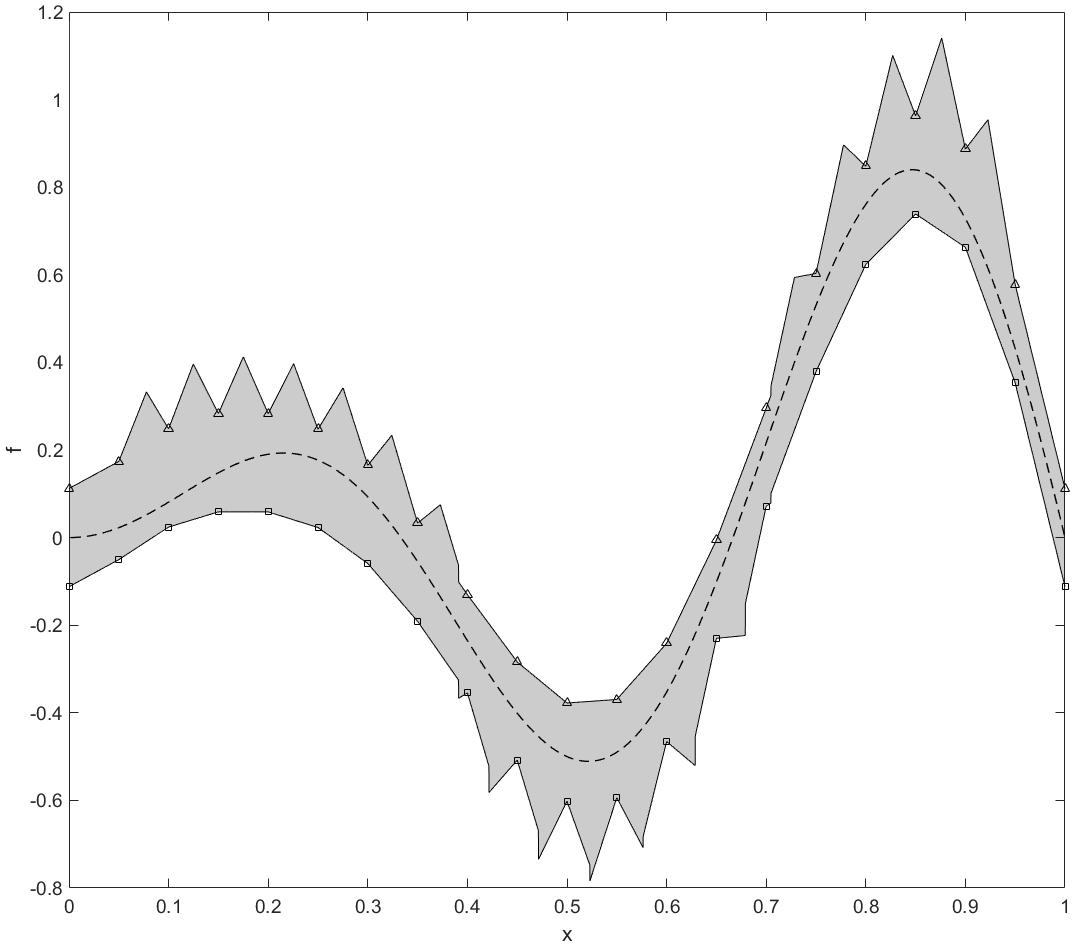}
\caption{ The lower $f^{low}(x)$ (line with squares)  and upper $f^{up}(x)$ (line with triangles) solutions; these can be compared with the exact source function $f^*(x)=x \sin{(3 \pi x)}$ (dashed line). }
\label{fig:convex2UPLOW}
\end{center}
\end{figure}

Finally, we uniformly select 21 random points from the reconstructed source function $f(x)$ (Fig. \ref{fig:sourcedeltaGAP10-3}), and using the formulas \eqref{convexLOW} and \eqref{convexUP} we construct the lower $f^{low}(x)$ and upper $f^{up}(x)$ solutions (see Fig. \ref{fig:convex2UPLOW}).

\section{Conclusions}
\label{Conclusion}

In this paper, by applying the asymptotic analysis, we propose a numerical-asymptotic approach to solving both forward and inverse problems of a nonlinear singularly perturbed PDE. The main advantage of this method is that it allows the approximation of the original high-order-differential-equation model with faster transiting internal layer through a simplified lower-order differential equation, which describes the solution of the problem over the entire domain of the problem definition except for a narrow region, the width of which is also estimated in this paper. This simplification will not decrease the accuracy of the inversion result, especially for inverse problems with noisy data, and thus provides a robust inversion solver~-- AER. We believe this approach can be applied to a wide class of asymptotically perturbed PDEs. The asymptotic analysis makes it possible to establish a simpler link relation between the input data and the quantity of interest in the inverse problems, which greatly simplifies the procedure for solving original PDE-based inverse problems.

\section{Acknowledgement}


This work has been supported by Beijing Natural Science Foundation (Key project No. Z210001), National Natural Science Foundation of China (No. 12171036), the Guangdong Fundamental and Applied Research Fund (No. 2019A1515110971) and Shenzhen National Science Foundation (No. 20200827173701001).

\bibliographystyle{model1-num-names}
\bibliography{sn-article}

\end{document}